\documentclass[11pt]{article}
\usepackage{amssymb}

\newtheorem{thm}{Theorem}[section]
\newtheorem{prop}[thm]{Proposition}
\newtheorem{rem}[thm]{Remark}
\newtheorem{lemma}[thm]{Lemma}
\newtheorem{cor}[thm]{Corollary}

\newtheorem{exam}[thm]{Example}
\newtheorem{ques}[thm]{Question}
\newtheorem{prob}[thm]{Problem}
\newtheorem{conj}[thm]{Conjecture}

\newcommand{\co}{\colon\thinspace}
\def\a{{\alpha}}
\def\b{{\beta}}

\def\Q{{\mathbb Q}}
\def\R{{\mathbb R}}
\def\Z{{\mathbb Z}}

\def\square{{\vcenter{\hrule height.4pt
       \hbox{\vrule width.4pt height5pt \hskip5pt
            \vrule width.4pt}
       \hrule height.4pt}}}

\def\qed{\hfill$\square$\bigskip}
\def\pf{{\noindent{\bf Proof \hspace{2mm}}}}

\begin{document}

\title{Orderable 3-manifold groups}
\author{S. Boyer\footnote{{Partially supported by grants: NSERC
OGP0009446 and FCAR ER-68657.}} ,
D. Rolfsen\footnote{{Partially supported by NSERC grant OGP0008082 and
Marsden fund UOA520.}} ,
B. Wiest\footnote{{Partially supported by NSERC grant OGP0008082 and a
PIMS postdoctoral
fellowship.}}}
\date{ }
\maketitle


\abstract{We investigate the orderability properties of fundamental groups 
of 3-dimensional manifolds.  Many 3-manifold groups support left-invariant
orderings, including all compact $P^2$-irreducible manifolds with positive
first Betti number.  For seven of the eight geometries (excluding 
hyperbolic) we are able to characterize which manifolds' groups support a
left-invariant or bi-invariant ordering.  We also show that manifolds 
modelled on these geometries
have virtually bi-orderable groups.  The question of virtual orderability
of 3-manifold groups in general, and even hyperbolic manifolds, remains
open, and is closely related to conjectures of Waldhausen and others.}


\section{Introduction}\label{S:Intro}

A group $G$ is called {\it left-orderable} (LO) if its elements can
be given a (strict) total ordering $<$ which is left invariant, meaning that
$g < h  \Rightarrow fg < fh$ if $f,g, h \in G$.  We will say that $G$ is
{\it bi-orderable} (O) if it admits a total ordering which is simultaneously
left and right invariant (historically, this has been called ``orderable'').
A group is called {\it virtually left-orderable}  or {\it virtually
bi-orderable} if it has a finite index subgroup with the appropriate property.

It has recently been realized that many of the groups which
arise naturally in topology are left-orderable.  Dehornoy provided a
left-ordering for the Artin braid groups
\cite{Dehornoy}; see also \cite{FGRRW} and \cite{SW}.  Rourke
and Wiest \cite{RW} extended this, showing that mapping class groups of all
Riemann surfaces with nonempty boundary (and possibly with punctures) are
left-orderable.  In general these groups are not bi-orderable.  On the other
hand, the pure Artin braid groups are known to be bi-orderable \cite{RZ},
\cite{KR}, and Gonzales-Meneses \cite{G-M} has constructed a
bi-ordering on the pure braid groups of orientable surfaces $PB_n(M^2)$.

The goal of the present paper is to investigate the orderability of
the fundamental groups of compact, connected $3$-manifolds, a class we
refer to as {\it $3$-manifold groups}.  We include nonorientable manifolds,
and manifolds with boundary in the analysis.  It will be seen that
left-orderability is a rather common property in this class, but is by no
means universal. After reviewing some general properties of orderable
groups in \S \ref{genord}, we begin our investigation of $3$-manifold
groups in \S \ref{gen3man}, asking not only if such a group is left- or
bi-orderable, but also if these properties hold virtually.
In other words, we examine whether or not there is a  finite cover of the
manifold whose group is left- or bi-orderable.  The following
is one of our general results.

\begin{thm} \label{hsext} 
Suppose that $M$ is a compact, connected and
$P^2$-irreducible $3$-manifold. \\
$(1)$ A necessary and sufficient condition that $\pi_1(M)$ be
left-orderable is that either $\pi_1(M)$ is trivial or there exists a
non-trivial homomorphism from $\pi_1(M)$ to a left-orderable group. \\
$(2)$ If $\pi_1(M)$ is not virtually left-orderable, then $M$ is closed,
orientable and homotopically atoroidal, that is, there is no $\mathbb 
Z \oplus \mathbb Z$ subgroup in $\pi_1(M)$.  
\end{thm}

Part (1) of this theorem follows from work of Howie and Short and an
observation of Boileau. See theorem \ref{homtoLO}.
Part (2) is a consequence of part (1) and work of Casson, Gabai, 
Jungreis, and Luecke. See the
discussion following conjecture \ref{vpos}.

Theorem \ref{hsext} implies that a compact, connected,
$P^2$-irreducible $3$-manifold $M$ whose first Betti number $b_1(M)$ is
larger than zero has a left-orderable fundamental group.
This is, in fact, the generic case, as it is well-known that $b_1(M) > 0$
when $M$ is neither a $3$-ball nor a $\Q$-homology $3$-sphere (cf. lemma
\ref{posb1}). On the other hand, we will see below that
not every $P^2$-irreducible
$\Q$-homology $3$-sphere $M$ has a left-orderable fundamental group,
even if this group is infinite.
Nevertheless, it frequently does up to taking a finite index
subgroup. Danny Calegari pointed out to us that this is the case when $M$
is atoroidal
and admits  a transversely orientable taut foliation owing to the existence
of a faithful
representation $\pi_1(M) \to Homeo_+(S^1)$ arising from Thurston's
universal circle
construction. More generally the following
result holds (see \S \ref{gen3man}).

\begin{prop} \label{hom3}
Let $M$ be an irreducible $\Q$-homology $3$-sphere and $\hat M \to M$
the finite-sheeted cover corresponding to the
commutator subgroup of $\pi_1(M)$. If there is a homomorphism 
$\pi_1(M) \to Homeo_+(S^1)$ with non-abelian image, then $\pi_1(\hat M)$ 
is left-orderable. In particular if $M$ is a $\Z$-homology $3$-sphere 
admitting such a homomorphism $\pi_1(M)\to Homeo_+(S)$, then $\pi_1(M)$ 
is left-orderable.
\end{prop}

Background material on Seifert fibred spaces is
presented in \S \ref{sfs} while in \S \ref{cod1} we examine the connection
between
orderability and codimension $1$ objects such as foliations. After these
general results, we focus attention on the special class of Seifert fibred
$3$-manifolds, possibly non-orientable, a convenient class which is
well-understood,
yet rich in structure. For this case we are able to supply complete answers.

\begin{thm} \label{SeifertLO}
For the fundamental group of a compact, connected, Seifert fibred
space $M$ to be left-orderable it is necessary and sufficient that $M 
\cong S^3$ or one of the following two sets of conditions holds: \\
$(1)$ $b_1(M) > 0$ and
$M \not \cong {P^2} \times S^1$; \\
$(2)$ $b_1(M) = 0$, $M$ is orientable, $\pi_1(M)$ is
infinite, the base orbifold of $M$
is of the form $S^2(\alpha_1, \alpha_2,$ $\ldots, \alpha_n)$,
and $M$ admits a horizontal foliation.
\end{thm}

The definition of horizontal foliation is given in \S \ref{cod1}. 
When applying this theorem, it is worth keeping in mind that Seifert 
manifolds whose first Betti number is zero and which have infinite 
fundamental group admit unique
Seifert structures (see \cite{Jc}, theorem VI.17). We also remark
that owing to the combined work of Eisenbud-Hirsch-Neumann \cite{EHN},
Jankins-Neumann \cite{JN2} and Naimi \cite{Na},
it is known exactly which Seifert manifolds admit horizontal foliations
(see theorem \ref{class}).
This work and theorem \ref{hsext} show that left-orderability is a
much weaker condition than the existence of a horizontal foliation
for Seifert manifolds of positive first Betti number.

Roberts and Stein have shown \cite{RS} that a necessary and sufficient
condition for an irreducible, non-Haken Seifert fibred manifold to
admit a horizontal foliation is that its fundamental group act
non-trivially (i.e. without a global fixed point)
on $\mathbb R$, a condition which is (in
this setting)  equivalent to the left-orderability of the group (theorem
\ref{linnell}).
Since these  Seifert manifolds have base orbifolds of the form
$S^2(\alpha_1, \alpha_2, \alpha_3)$,
theorem \ref{SeifertLO} can be seen as a generalization of their
result.

Theorem \ref{SeifertLO} characterizes the Seifert manifold groups which 
are left-orderable. In order to characterize those which are bi-orderable, 
we must first deal with the same question for surface groups. It is 
well-known that free groups are bi-orderable. Moreover, it was observed by
Baumslag that the fundamental group of an orientable surface
is residually free, and therefore bi-orderable (see \cite{Smythe} and
\cite{Baumslag}). In \S \ref{biosurf} we give a new proof
of the bi-orderability of closed orientable surface groups, and settle the
orderability question for closed, nonorientable surface groups.  This 
result also appears in \cite{RolW}.

\begin{thm}\label{bisur}
If $M$ is any connected surface other than the projective plane or 
Klein bottle, then $\pi_1(M)$ is bi-orderable.
\end{thm}

In \S \ref{bioseif} we will use this result to prove

\begin{thm}\label{SeifertBO}
For the fundamental group of a compact, connected, Seifert fibred
space $M$ to be bi-orderable, it is necessary and sufficient that either \\
$(1)$ $M$ be homeomorphic to one of $S^3, S^1 \times S^2, S^1 \tilde \times
S^2$ $($the non-trivial $2$-sphere bundle over the circle$)$, a solid Klein
bottle, or \\
$(2)$ $M$ be the total space of a locally trivial, orientable circle bundle
over a surface other than $S^2, P^2$, or the Klein bottle.
\end{thm}

\begin{cor}
The fundamental group of any compact Seifert fibred manifold is virtually
bi-orderable. \end{cor}

\pf A Seifert manifold is always finitely covered by an orientable Seifert
manifold with no exceptional fibres, that is, a locally trivial
circle bundle over an orientable surface. If that surface happens to be 
a $2$-sphere, there is a further finite cover whose total space is either 
$S^3$ or $S^1 \times S^2$.
\qed

Seifert manifolds account for six of the eight $3$-dimensional geometries.
Of the two remaining geometries, hyperbolic and $Sol$, the latter
is fairly simple to understand in terms of orderability properties.  
In \S \ref{osol} we prove the following theorem.

\begin{thm} \label{Solchar}
Let $M$ be a compact, connected $Sol$ manifold. Then  \\
$(1)$ $\pi_1(M)$ is left-orderable if and only if $\partial M \ne
\emptyset$, or
$M$ is non-orientable, or $M$ is orientable and a torus bundle over the
circle. \\
$(2)$ $\pi_1(M)$ is bi-orderable if and only if $\partial M \ne \emptyset$ 
but $M$ is not a twisted $I$-bundle over the Klein bottle, or $M$
is a torus bundle over the circle whose monodromy in $GL_2(\mathbb Z)$ has
at least one positive eigenvalue. \\
$(3)$ $\pi_1(M)$ is virtually bi-orderable.
\end{thm}

In a final section we consider hyperbolic 3-manifolds.  This is the 
geometry in which the orderability question seems to us to be the most 
difficult, and we have only partial results.  We discuss a very recent 
example \cite{RSS} of a closed hyperbolic 3-manifold whose fundamental 
group is not left-orderable.
On the other hand, there are many closed hyperbolic 3-manifolds whose 
groups {\it are} LO -- for example those which have infinite first
homology (by theorem \ref{hsext}).
This enables us to prove the following result.

\begin{thm} \label{geoms}  For each of the eight 3-dimensional geometries,
there exist closed, connected, orientable $3$-manifolds with the given
geometric
structure whose fundamental groups are left-orderable. There are also
closed, connected, orientable $3$-manifolds with the given geometric
structure whose groups are not left-orderable.
\end{thm}

This result seems to imply that geometric structure and orderability are
not closely related. Nevertheless compact, connected hyperbolic 
$3$-manifolds are conjectured to have finite covers with positive first 
Betti numbers, and if this is true, their fundamental groups are virtually
left-orderable (cf. corollary \ref{gencond} (1)). One can also ask whether 
they have finite covers with bi-orderable groups, though to put the 
relative difficulty of this question in perspective, note that nontrivial, 
finitely generated, bi-orderable groups have
positive first Betti numbers (cf. theorem \ref{OLI}).

We close the introduction by listing several questions and problems arising from this study.

\begin{ques}
{\rm Is the fundamental group of a compact, connected, orientable 
$3$-manifold virtually left-orderable?  What if the manifold is hyperbolic?}
\end{ques}

It is straightforward to argue that $3$-manifold groups are virtually
torsion free -- the main ingredient in the proof is the fact that prime
orientable 3-manifolds with torsion in the fundamental group have finite
fundamental group. \\

We saw in theorems \ref{SeifertBO} and \ref{Solchar} that the
bi-orderability of the fundamental groups
of Seifert manifolds and Sol manifolds  can be detected in a
straightforward manner.
The same problem for hyperbolic manifolds appears to be much more subtle.

\begin{ques}
{\rm Is there a compact, connected, orientable irreducible $3$-manifold
whose fundamental group is not
virtually bi-orderable?  What if the manifold is hyperbolic?}
\end{ques}

\begin{prob}
{\rm Find necessary and sufficient conditions for the fundamental group of
a compact, connected $3$-manifold which fibres over the circle to be
bi-orderable. Equivalently, can one find bi-orderings of free groups or
surface groups which are invariant under the automorphism corresponding
to the monodromy of the fibration?}
\end{prob}

This problem is quickly dealt with in the case when the fibre is of 
non-negative
Euler characteristic,
so the interesting case involves fibres which are hyperbolic surfaces. When
the boundary of the surface is
non-empty, Perron and Rolfsen \cite{PR} have found a sufficient condition for
bi-orderability; for instance, the fundamental group of the figure eight
knot exterior has a bi-orderable fundamental group.\\

\noindent {\bf Acknowledgements } This paper owes much to very useful
discussions with Hamish Short and Michel Boileau.  We also thank 
Gilbert Levitt, Denis Sjerve, and Danny Calegari for helpful 
comments. The referee made a number of suggestions which greatly 
improved the exposition.


\section{Ordered and bi-ordered groups.} \label{genord}

We summarize a few
facts about left-orderable (LO) groups, bi-orderable (O) groups  and other
algebraic matters.  Good references are \cite{MR}, \cite{Passman}.

\begin{prop}
If $G$ is left-orderable, then $G$ is torsion-free.
\end{prop}

\pf  If $g \ne 1$, we wish to show $g^p \ne 1$.  Without loss of 
generality, $1 < g$.  Then $g < g^2$, $g^2 < g^3$, etc. so an easy 
induction shows $1 < g^p$ for all positive $p$. \qed

Thus we can think of left-orderability as a strong form of torsion-freeness.
The following lemma will be crucial for our classification of Seifert-fibred
spaces with bi-orderable fundamental group.

\begin{lemma} \label{central}
In a bi-orderable group $G$, a non-zero power of an element $\gamma$ is
central if and only if
$\gamma$ is central.
\end{lemma}

\pf  Obviously, in any group, powers of a central element are central.
On the other hand, suppose there is an integer $n > 0$ such that
$\gamma^{n}$ is central in $G$. If there is some $\mu \in G$
which does not commute with $\gamma$, say $\gamma \mu \gamma^{-1} <  \mu$.
Then by invariance under conjugation,
$\gamma^2 \mu \gamma^{-2} < \gamma \mu \gamma^{-1} < \mu$ and by induction
$\gamma^k \mu \gamma^{-k} < \mu$ for each positive integer
$k$.  We arrive at the contradiction:
$\mu = \gamma^{n} \mu \gamma^{-n} < \mu$.  The case
$\gamma \mu \gamma^{-1} >  \mu$ similarly leads to a contradiction.
Hence $\gamma$ must be central in $\pi_1(M)$.
\qed

A group $G$ is $LO$ if and only if there exists a subset $P \subset G$ (the
positive cone) such that (1) $P\cdot P = P$ and (2) for every $g\ne 1$
in $P$, exactly one of $g$ or $g^{-1}$ belongs to $P$.  Given such a $P$, 
the recipe $g<h$ if and only $g^{-1}h \in P$ is easily seen to define a
left-invariant strict total order, and conversely such an ordering defines 
the set $P$ as the set of elements greater than the identity.  The group 
$G$ is bi-orderable if and only if it admits a subset $P$ satisfying 
(1), (2), and in addition (3) $gPg^{-1} \subset P$ for all $g \in G$.

As we shall see in a moment, the class of $LO$ groups is closed under 
taking subgroups, extensions, directed unions, and free products. The 
class of $O$ groups is also invariant under taking subgroups, directed 
unions and free products, but not necessarily under extensions. An 
instructive example is the fundamental group of the Klein bottle:
$$G = \langle m,l\ ; \  lml^{-1} = m^{-1} \rangle.$$
This contains a normal subgroup $\Z$ generated by $m$, and the quotient 
$G/\Z$ is also an infinite cyclic group.  Of course $\Z$ is bi-orderable, 
so the extension $G$ of $\Z$ by $\Z$ is certainly left-orderable, by 
lemma \ref{oext} below. However, $G$ is not biorderable,
for if we had a biorder with $m>1$ then it would follow that
$1<lml^{-1}=m^{-1}<1$; if $m<1$ a similar contradiction arises.

\begin{lemma}[Orderability of extensions] \label{oext}
Let $f\co G \to H$ be a surjective homomorphism of groups, with kernel $K$,
and assume that both $H$ and $K$ are left-ordered, with positive cones $P_H$
and $P_K$, respectively.  Then the subset $P = P_K \cup f^{-1}(P_H)$
defines a left-invariant ordering on $G$ by the rule
$g<g' \Leftrightarrow g^{-1}g' \in P$. If $H$ and $K$ are bi-ordered, and
if in addition $P_K$ is normal in $G$, then this rule defines a bi-ordering
of $G$. \end{lemma}

\pf Routine, and left to the reader.
\qed

A {\it left action} of a group $G$ on a set $X$ is a homomorphism
$\phi$ from $G$ to the permutation group of $X$. For $g \in G$ and 
$x \in X$ we denote $\phi(g)(x)$ by $g(x)$. If $1 \in G$ is the
only group element that acts as the identity on $X$, the action is said 
to be {\it effective}.

\begin{thm} {\rm (Conrad, 1959)} \label{eff}
A group $G$ is left-orderable if and only if it acts effectively on a
linearly ordered set $X$ by order-preserving bijections.
\end{thm}

\pf  One direction is obvious, as a left-ordered group acts upon itself via
multiplication on the left.  On the other hand, assume $G$ acts on $X$ in such
a way that for every $g \in G$, $x<y \Leftrightarrow g(x)<g(y)$.
Let $\prec$ be some {\it well}-ordering of the elements of $X$,
completely independent of the given ordering $<$ and of the $G$-action
(such an order exists, by the axiom of choice).
Compare $g \ne h \in G$ by letting $x_0 \in X$ be the smallest $x$, in the
well-ordering $\prec$, such that $g(x) \ne h(x)$.  Then say that $g<h$ or
$h<g$ according as $g(x_0)<h(x_0)$ or $h(x_0)<g(x_0).$  This can easily be
seen to be a left-invariant ordering of $G$. \qed

Thus, the group $Homeo_+(\R)$ of order-preserving bijections is
left-orderable; it acts effectively on $\R$ by definition. It follows that
the universal covering group $\widetilde{SL}_2(\R)$ of $PSL_2(\R)$ is
left-orderable,
a fact first noted by Bergman \cite{bergman}, as it acts effectively and
order-preservingly on
the real line $\R$.

Next we state a classical result.
A left-ordering of the group $G$ is said to be $Archimedian$
if for each  $a, b \in G$ with $1 < a < b$, there is a positive integer $n$
such that  $b < a^n$.

\begin{thm}  {\rm (Conrad 1959, H\"older 1902)}  If a left-ordered group
$G$ is Archimedian, then the ordering is necessarily bi-invariant.  Moreover
$G$ is isomorphic (in both the algebraic and order sense) with a subset of
the additive real numbers, with the usual ordering.  In particular, $G$ is
abelian. \qed
\end{thm}

This result simply implies that most interesting left-ordered groups are
non-archimedian.  The following offer alternative criteria for
left-orderability; this is well-known to experts -- see \cite{Linnell} for 
a proof.

\begin{thm}\label{linnell}
If $G$ is a countable group, then the following are equivalent: \\
$(1)$ $G$ is left-orderable. \\
$(2)$ $G$ is isomorphic with a subgroup of $Homeo_+(\R)$. \\
$(3)$ $G$ is isomorphic with a subgroup of $Homeo_+(\Q)$.
\qed
\end{thm}

It is known that orderability is a local property.  That is, a group is LO 
(resp O) if and only if every finitely generated subgroup is LO (resp O).
Closely related to this, we have

\begin{thm} {\rm (Burns-Hale \cite{BH})} \label{BH}
A group is left-orderable if and only if every non-trivial finitely 
generated subgroup has a non-trivial quotient which is left-orderable.
\qed
\end{thm}

We recall a definition due to Higman: a group is {\it locally-indicable}
(LI) if every nontrivial finitely-generated subgroup has $\Z$ as a quotient.
The following is also well-known to experts \cite{Co}, \cite{BH}, \cite{MR}.

\begin{thm} \label{OLI}
If $G$ is a bi-orderable group, then $G$ is locally indicable.
If $G$ is locally indicable, then $G$ is left-orderable.
Neither of these implications is reversible.
\end{thm}

\pf  If $G$ is bi-ordered, consider a finitely generated subgroup
$H = \langle h_1, \dots, h_k \rangle$, with notation chosen so that
$1 < h_1 < \dots < h_k$.  We recall that a subgroup $C$ is called 
{\it convex} if $f<g<h, f\in C, h\in C \Rightarrow g \in C$.  The convex 
subgroups of a left-ordered group are totally ordered by inclusion and 
closed under intersections and
unions.   Now, considering $H$ itself as a finitely generated left-ordered
group, we let $K$ be the union of all convex subgroups of $H$ which do not
contain $h_k$.  Then one can use bi-orderability and a generalization of the
Conrad-H\"older theorem (or see \cite{Co} for a more general argument) to show that $K$ is normal in $H$, and the quotient $H/K$ is isomorphic with 
a subgroup of
$(\R, +)$.  Being finitely generated, $H/K$ is therefore isomorphic with a
sum of infinite cyclic groups, and so there is a nontrivial
homomorphism $H \to H/K \to \Z$, completing the first half of the theorem.

The second half follows directly from theorem \ref{BH}, and the
observation that $\Z$ is left-orderable. Finally, the fact that neither
implication is reversible is discussed in the paragraph which follows.
  \qed

Bergman \cite{bergman} observed that even though $\widetilde{SL}_2(\R)$ is
left-orderable, it is not locally-indicable: for example, it contains the 
perfect infinite group $\langle x, y, z: x^2 = y^3 = z^7 = xyz \rangle$, 
which happens to be the fundamental
group of a well-known homology sphere. The braid groups $B_n$, for
$n > 4$ are further examples of LO groups which are not locally indicable, as
their commutator subgroups $[B_n,B_n]$ are finitely generated and perfect
(see \cite{GorinLin}).  The
braid groups $B_3$ and $B_4$, and the Klein bottle group provide examples of
locally-indicable groups which are not bi-orderable.  There is a
characterization
of those left-orderable groups which are locally indicable in \cite{RR}. For
instance for solvable groups \cite{ChiswellKropholler}, and more generally,
elementary amenable groups \cite{Linnell}, the concepts of
left-orderability and local indicability coincide.

Our analysis of the orderability of the 
fundamental groups of compact $3$-manifolds will also rely heavily on the 
next two results. 

\begin{prop} {\rm (Vinogradov \cite{vinogradov})} 
\label{vin}
A neccesary and sufficient condition for a free product $G = G_1 * G_2 *
\ldots * G_n$ of
groups to be left-orderable, respectively bi-orderable, is that each $G_j$
has this
property.
\qed
\end{prop}

\begin{prop} \label{virtvin}
A neccesary and sufficient condition for a free product $G = G_1 * G_2 *
\ldots * G_n$ of
groups to be  virtually left-orderable (resp. virtually bi-orderable) is
that each $G_j$ have this
property.
\end{prop}

\pf For each $j$, the intersection of a finite index LO subgroup of 
$G$ with $G_j$ is a finite index LO subgroup of $G_j$. Thus $G_j$ is 
virtually left-orderable.

On the other hand if each $G_j$ is virtually left-orderable, there are
surjective homomorphisms $\phi_j\co G_j \to F_j$ where
$F_1, F_2, \ldots , F_n$ are finite groups and $\mbox{ker}(\phi_j)$ is
left-orderable. By the Kurosh subgroup theorem \cite{ScWa}, the kernel of
the obvious homomorphism
$G_1 * G_2 * \ldots * G_n \to F_1 \times F_2 \times \ldots \times F_n$ is a
free product of a free group and groups isomorphic to
$\mbox{ker}(\phi_1), \ldots,\mbox{ker}(\phi_n)$.
This finite-index subgroup is left-orderable by proposition \ref{vin}.

A similar argument shows the analoguous statement for bi-orderable
groups. \qed

In fact, the previous two results also hold for free products of 
infinitely many groups. This follows from the fact that orderability 
and bi-orderability are local properties, together with the Kurosh 
subgroup theorem.

We mention in passing the following theorem of Farrell, which relates 
orderability with covering space theory.

\begin{thm} {\rm (Farrell \cite{Fa})}
Suppose $X$ is a locally-compact, paracompact topological space, and let
$p\co \widetilde{X} \to X$ the universal covering.  Then
$\pi_1(X)$ is left-orderable if and only if there is a topological embedding
$h\co\tilde{X} \to X \times \R$ such that $pr_X h = p$.
\qed
\end{thm}

We conclude this section with certain facts about orderable groups, which
makes orderability properties worthwhile knowing.  Of particular interest are
the deep properties of the group ring $\Z G$.

\begin{thm} {\rm (see eg. \cite{Passman})} If $G$ is left-orderable,
then $\Z G$ has no zero divisors, and
only the units $ng$ where $n$ is a unit of $\Z$ and $g \in G$.  The same 
is true for any integral domain $R$ replacing $\Z$.
\qed
\end{thm}

A proof is not difficult, the idea being to show that in a formal product, the largest (and smallest) terms in the product cannot be cancelled by 
any other term.  The conclusions of this theorem are conjectured to
be true for arbitrary torsion-free groups.  For bi-orderable groups we know
even more.

\begin{thm} {\rm (Mal'cev \cite{malcev}, B. Neumann \cite{neumann})}
If $G$ is bi-orderable then $\Z G$ embeds in a division algebra.
\qed
\end{thm}

\begin{thm} {\rm (LaGrange, Rhemtulla \cite{LR})}
Suppose $G$ and $H$ are groups with $G$ left-orderable.  Then $G$ and $H$ are
isomorphic groups if and only if their group rings $\Z G$ and $\Z H$ are
isomorphic as rings.
\qed
\end{thm}


\section{General remarks on ordering 3-manifold groups} \label{gen3man}

\subsection{Orderability}

Good references for the basic facts on  $3$-manifolds that we shall use 
in this paper are \cite{He} and \cite{Jc}.

Recall that a compact, connected $3$-manifold $M \ne S^3$ splits into a
product of prime $3$-manifolds under connected sum
$$M \cong M_1 \# M_2 \# \ldots \# M_n$$
(see theorem 3.15 of \cite{He}).
Clearly then $\pi_1(M) \cong  \pi_1(M_1) * \pi_1(M_2) * \ldots * \pi_1(M_n)$ 
is LO, respectively O, if and only if each $\pi_1(M_j)$ has this property
(cf. proposition \ref{vin}).
Since the fundamental group of a prime, reducible $3$-manifold is 
$\mathbb Z$, it suffices to investigate the orderability of the groups 
of irreducible $3$-manifolds. We can specialize further. 
Recall that a 3-manifold is $P^2$-irreducible if and
only if it is irreducible and contains no $2$-sided $P^2$. For an 
irreducible manifold, containing a $2$-sided $P^2$ is equivalent to 
the manifold being non-orientable and having a $\Z/2$ subgroup in its 
fundamental group (\cite{Ep}, Theorem 8.2).  Therefore we need only 
consider $P^2$-irreducible 3-manifolds.    

The method of proof of lemma 2 of  \cite{HoSh} yields the following 
result.

\begin{thm}\label{homtoLI} {{\rm (Howie-Short)}}
Suppose that $M$ is a compact, connected, $P^2$-irreducible 3-manifold 
and that $\pi_1(M)$ is nontrivial. A necessary and sufficient condition 
that $\pi_1(M)$ be locally indicable is that $b_1(M) > 0$. 
\qed
\end{thm}
 
An analogous result holds in the situation of interest 
to us (cf. theorem \ref{hsext}(1)).

\begin{thm}\label{homtoLO}
Suppose that $M$ is a compact, connected, $P^2$-irreducible 3-manifold and
that $\pi_1(M)$ is nontrivial. A necessary
and sufficient condition that $\pi_1(M)$ be left-orderable is that there
exists a homomorphism $h\co \pi_1(M) \to L$  onto a nontrivial 
left-orderable group $L$.
\end{thm}

Before proving theorem \ref{homtoLO} we need the following
well-known lemma.

\begin{lemma} \label{posb1}
If $M$ is a compact 3-manifold and either $M$ is closed and 
non-orientable, or $\partial M$ is nonempty but contains no
$S^2$ or $P^2$ components, then $b_1(M) > 0$. In particular this 
holds for nonorientable $P^2$-irreducible 3-manifolds.
\end{lemma}

\pf  We wish to show that $b_1(M)$
is at least $1$.  Noting that closed 3-manifolds have
zero Euler characteristic, if $W$ is the double of $M$, then 
$0 = \chi(W) = 2 \chi(M) - \chi(\partial M)$, and
so $2 \chi(M) = \chi(\partial M)$. Our hypotheses imply that $H_3(M) = 0$
while each component of $\partial M$ has a non-positive Euler
characteristic, thus
$$0 \geqslant {1 \over 2}\chi(\partial M) = \chi(M) = \sum (-1)^j b_j(M) = 
1 - b_1(M) + b_2(M) $$  and we conclude that
$b_1(M) \geqslant1 + b_2(M) \geqslant1.$
\qed

\noindent {\bf Proof of theorem \ref{homtoLO}} $\;$
Necessity is obvious.  For sufficiency, assume there is a
surjection   $h\co \pi_1(M) \to L$, with $L$ nontrivial left-orderable. We
wish to show that $\pi_1(M)$ is left-orderable. Using the Burns-Hale
characterization (theorem \ref{BH}), it suffices to show that every
nontrivial finitely generated subgroup $H$ of $\pi_1(M)$ has a homomorphism
onto a nontrivial left-ordered group.  Consider such a group $H$
and distinguish two cases.  If $H$ has finite index in $\pi_1(M)$, then
$h(H)$ is a finite index subgroup of $L$ and therefore nontrivial.  So in
this case we can just take the restriction of $h$ to $H$.

Now suppose $H$ has infinite index and let $p\co \tilde{M} \to M$ be the
corresponding covering space, i.e.\ $p_\#(\pi_1(\tilde{M},\tilde{*})) = H$.
Although $\tilde{M}$ is necessarily noncompact, by a theorem of Scott 
\cite{Sc1}, there is a compact submanifold $C \subset \tilde{M}$ whose 
fundamental group is isomorphic, via inclusion, with $\pi_1(\tilde{M}).$  
The manifold $C$ must have
nonempty boundary, otherwise it would be all of $\tilde{M}$.  Suppose
that $S \subset \partial C$ is a 2-sphere. Since $M$ is irreducible,
so is $\tilde{M}$ (see \cite{MSY, Dunwoody, Ha}), and therefore $S$ bounds
a 3-ball $B$ in $\tilde{M}$. We claim that $B \cap C = S$, for otherwise
we would have $C \subset B \subset \tilde{M}$, contradicting that the
inclusion of $C$ in $\tilde{M}$ induces an isomorphism of nontrivial 
groups.  Thus we may attach $B$ to $C$ without
changing the property that $i_*\co \pi_1(C) \to \pi_1(\tilde{M})$ is an
isomorphism.  Therefore we may assume $\partial C$ contains no 2-sphere
components.  Next we wish to show that no component of $\partial C$ is
a projective plane.   If there were such a component, it would contain a
loop $\alpha$ which reverses orientation of $\tilde{M}$, and hence is
nontrivial in $\pi_1(\tilde{M}).$  On the other hand, since it lies in the
projective plane, $\alpha^2$ = 1; which would imply that $\pi_1(\tilde M)$,
and therefore $\pi_1(M)$ has an
element of order 2, which is not allowed.   We now have that $\partial C$ is
nonempty, but contains no spheres or projective planes.  By lemma
\ref{posb1}, $H_1(C)$ is infinite, and therefore maps onto $\Z$.  Preceding 
this homomorphism by the Hurewicz map $\pi_1(C) \to H_1(C)$ gives the 
required homomorphism of $H$ onto $\Z$. \qed

\begin{cor} \label{gencond}
Let $M$ be a compact, connected, prime $3$-manifold, possibly with boundary. 
\\
$(1)$ If $M$ is orientable with $b_1(M) > 0$, then $\pi_1(M)$ is
left-orderable.
\\
$(2)$ If $M$ is non-orientable then $\pi_1(M)$ is left-orderable if and
only if it contains no element of order $2$.
\end{cor}

\pf If $M$ is reducible, its group is $\mathbb Z$, so the corollary holds.
On the other hand if it is irreducible the result follows from  theorem
\ref{homtoLO} and lemma \ref{posb1}.
\qed

\begin{cor}\label{knotLO}
Let $G = \pi_1(S^3 \setminus K)$ be a knot or link group.  Then $G$
is left-orderable.
\end{cor}

\pf The only point to observe is that the group of a split link is a 
free product of
the groups of non-split links, whose complements are
irreducible (cf. lemma \ref{vin}). 
\qed

It follows from theorem \ref{homtoLI}, lemma \ref{posb1}, and corollary \ref{gencond} that the only compact prime 3-manifolds which
can have LO but not LI fundamental groups are those with finite first
homology. Bergman's abovementioned example -- a homology sphere whose
fundamental group is contained in $\widetilde{SL}_2(\R)$ -- is just such a manifold.

We saw above that compact, connected, orientable, irreducible $3$-manifolds
with positive first Betti
numbers have left-orderable groups. Such manifolds are Haken. On the other
hand, not all
Haken $3$-manifolds have left-orderable groups (see eg. theorem
\ref{SeifertLO}). The simplest examples
were constructed by Boileau, Short and Wiest.

\begin{exam}
{\rm (Boileau, Short and Wiest) Let $X$ be the exterior of a trefoil 
knot $K \subset S^3$ and let $\mu, \phi$ denote, respectively, the 
meridional slope on $\partial X$ and the slope corresponding to a fibre 
of the Seifert structure on $X$. Fix a base point$* \in \partial X$
and oriented representative curves $C_1, C_2$ for $\mu$ and $\phi$ based 
at $*$. The group $\pi_1(X; *)$ has a presentation
$\langle  x,y \; | \; x^2 = y^3 \rangle$ where $xy^{-1}$ represents the
class of $C_1$ while
$x^2$ represents that of $C_2$. Since $C_1$ and $C_2$ intersect once
algebraically,
there is a homeomorphism $f\co \partial X \to \partial X$ which switches
them. The manifold
$M = X \cup_f X$ is Haken, because the separating torus is incompressible.
We claim that its fundamental group is not left-orderable.

Assume to the contrary that $<$ is a left-order on
$$\pi_1(M;*) = \langle x_1, y_1, x_2, y_2 \; | \; x_1^2 = y_1^3, x_2^2 =
y_2^3, x_1^2 = x_2y_2^{-1},
x_2^2 = x_1y_1^{-1} \rangle.$$
Without loss of generality, $x_1 > 1$. The relation $y_1^3 = x_1^2$ implies
that $y_1 > 1$ as well.
Hence $x_1 > y_1^{-1}$, or equivalently, $x_1^ 2 > x_1y_1^{-1}$. If $x_2 >
1$, a similar argument shows
$x_2^ 2 > x_2y_2^{-1}$. But then $x_1^ 2 > x_1y_1^{-1} = x_2^ 2 >
x_2y_2^{-1} = x_1^2$, a contradiction.
Hence $x_2 < 1$. Now $y_1^2 > x_1^{-1}$ implies 
$x_1^2 = y_1^3 > y_1x_1^{-1}$,
while  $x_2 < 1$ implies $x_2^{2} = y_2^3 < y_2 
< y_2x_2^{-1}$. Since $x_2^2$ commutes with $y_2 x_2^{-1}$, we deduce 
$x_2^{-2} > x_2 y_2^{-1} = x_1^2$. But then, 
$x_1^ 2 > y_1x_1^{-1} = x_2^ {-2} > x_1^2$, 
another contradiction.
It follows that there is no left-order on $\pi_1(M)$. }
\end{exam}


\subsection{An application to mappings between 3-manifolds}

Now that we have an example of a 3-manifold whose group is infinite and
torsion-free, yet not left-orderable (there are many others), it is 
appropriate to point out an application of theorem \ref{hsext}. 
An important question in 3-manifold theory is whether, given two closed 
oriented 3-manifolds $M$ and $N$, there exists a degree one map $M \to N$,
or, more generally, a map of nonzero degree.  The following can be viewed 
as providing a new ``obstruction'' to the existence of such a map.

\begin{thm} \label{degree}
Let $M$ and $N$ be closed, oriented 3-manifolds, with $M$ prime.  
Suppose $\pi_1(N)$ is nontrivial and left orderable, but $\pi_1(M)$ is 
\emph{not} left orderable. Then any mapping $M \to N$ has degree zero.  
\end{thm}

\pf Being prime and orientable, $M$ is either irreducible
or $S^2 \times S^1$, but the latter possibility is excluded by hypothesis.
Suppose there were a mapping $M \to N$ of nonzero degree.  According 
to the lemma below, the induced map $\pi_1(M) \to \pi_1(N)$ would be 
nontrivial.  But then, by theorem \ref{homtoLO}, $\pi_1(M)$ would be 
left-orderable, a contradiction. \qed

\begin{lemma}
If $f:M \to N$ is a mapping of nonzero degree, then $f_*(\pi_1(M))$ has
finite index in $\pi_1(N).$
\end{lemma}

\pf  Let $p: \tilde{N} \to N$ denote the cover corresponding to
$f_*(\pi_1(M))$, so there is a lift $\tilde{f}:M \to \tilde{N}$.  Now $\tilde{N}$ must be compact, otherwise $H_3(\tilde{N}) = 0$, and since 
$f$ factors through $\tilde{N}$ its degree would be zero. Thus the covering 
is finite-sheeted, and the index is finite. \qed


\subsection{Virtual orderability}
Though a $3$-manifold group may not be left-orderable, it seems likely that it contains a finite index subgroup which is.
We consider then, the virtual orderability properties of the group of a
prime $3$-manifold $M$. It is clear that we may restrict our attention to 
prime $3$-manifolds which are irreducible. Recall the following variant of 
a conjecture of Waldhausen.

\begin{conj} \label{vpos}
{\rm If $M$ is a compact, connected, $P^2$-irreducible $3$-manifold with
infinite fundamental group, then there is a finite cover $\tilde M \to M$ 
with $b_1(\tilde M) > 0$.}
\end{conj}

\noindent If this conjecture turns out to be true, then theorem \ref{homtoLO}
implies that any prime $P^2$-irreducible $3$-manifold $M$ has a virtually
left-orderable group. While examining the virtual
left-orderability of $\pi_1(M)$, we may as well assume that $M$ is 
closed and orientable
(corollary \ref{gencond}).  Under these conditions, if $\pi_1(M)$ 
contains a $\mathbb Z \oplus \mathbb Z$ subgroup, then $M$ is either 
Seifert fibred or admits a $\pi_1$-injective torus (\cite{CJ}, 
\cite{Ga2}). It is well-known that Conjecture \ref{vpos} holds in the 
former case while John Luecke has shown \cite{Lu} that it holds in 
the latter. In either case, $\pi_1(M)$ is virtually
left-orderable.   We may therefore assume that $M$ is homotopically
atoroidal as well as being irreducible and closed. 
Such manifolds are simple (they contain no essential surfaces of 
non-negative Euler characteristic) and conjecturally hyperbolic.  

We will see in \S \ref{bioseif} and \S \ref{osol} that the groups of 
Seifert and Sol
manifolds are virtually bi-orderable, but we do not know if this holds 
for hyperbolic manifolds. Remark that by theorem \ref{OLI}, if $M$ is a 
compact, connected,
$P^2$-irreducible $3$-manifold which has a virtually bi-orderable 
fundamental group, then this group is virtually locally
indicable. Hence it has a virtually positive first Betti number.  This
puts the relative difficulty of the virtual bi-orderability of 
$3$-manifold groups in perspective.

Next we apply theorem \ref{homtoLO} to prove proposition
\ref{hom3}. We begin with a simple lemma pointed out to us by Danny Calegari.

\begin{lemma}
Let $\Gamma$ be a group such that $H_2(\Gamma) = 0$. Suppose that
$1 \to A \to \tilde G \to G \to 1$ is a central extension of a group $G$ by
a group $A$. If
$\rho: \Gamma \to G$ is a homomorphism, then the restriction $\rho|[\Gamma,
\Gamma]$
lifts to a homomorphism $[\Gamma, \Gamma] \to \tilde G$.
\end{lemma}

\pf Fix $\gamma = \Pi_{i=1}^n [\alpha_i, \beta_i] \in [\Gamma, \Gamma]$
and let $\tilde a_i, \tilde b_i$ be arbitrary lifts of $\rho(\alpha_i),
\rho(\beta_i)$ to $\tilde G$. The centrality of the extension shows that
$\tilde \gamma := \Pi_{i=1}^n [\tilde a_i, \tilde b_i]$ is independent of
our choice of lifts. Now we claim that $\tilde \gamma$ is independent 
of the way we expressed $\gamma$ as a product of commutators. 
Equivalently, we claim that if $\Pi_{i=1}^n [\alpha_i, \beta_i] = 1$, 
then $\Pi_{i=1}^n [\tilde a_i, \tilde b_i] = 1 \in \tilde G$. Once we 
show this, the correspondence $\gamma \mapsto \tilde \gamma$ provides 
the desired lift of $\rho|[\Gamma, \Gamma]$.

Let $1 \to R \to F \stackrel{\phi}{\to} \Gamma \to 1$ be a free
presentation of $\Gamma$ and fix a lift $\tilde \phi: F \to \tilde G$ of
$\rho \circ \phi$. Then $\tilde\phi(R) \subset A$ lies in the centre of
$\tilde G$. Choose $x_i, y_i \in F$ which are sent to $\alpha_i, \beta_i$ by
$\phi$.  Then by construction, $\Pi_{i=1}^n [x_i, y_i] \in R \cap [F,F]$. On
the  other hand Hopf's formula [HiSt] shows that
$0 = H_2(\Gamma) = (R \cap [F,F])/[F,R]$, and so
$\Pi_{i=1}^n [x_i, y_i] = \Pi_{j=1}^m [f_j, r_j]$ for some
$f_j \in F$ and $r_j \in R$. Then
$$\Pi_{i=1}^n [\tilde a_i, \tilde b_i] =
\Pi_{i=1}^n [\tilde \phi(x_i), \tilde \phi(y_i)] =
\Pi_{j=1}^m [\tilde \phi(f_j), \tilde \phi(r_j)] = 1$$
since $\tilde \phi(r_j)$ is contained in the centre of $\tilde G$ for each
$j$. This completes the proof.
\qed

\noindent{\bf Proof of proposition \ref{hom3} }
Recall that $M$ is an irreducible $\Q$-homology
$3$-sphere and that $\hat M \to M$ is the cover corresponding to the
commutator subgroup of $\pi_1(M)$. We are
given a homomorphism $\rho\co \pi_1(M) \to Homeo_+(S^1)$ whose image is not
abelian and we want to deduce that $\pi_1(\hat M)$ is left-orderable.

Consider the central $\Z$ extension
$$ 1 \to \Z \to \widetilde{Homeo}_+(S^1) \to Homeo_+(S^1) \to 1$$
where $\widetilde{Homeo}_+(S^1)$ is the universal covering group of
$Homeo_+(S^1)$. This covering group can be identified
with the subgroup of $Homeo_+(\R)$ consisting of homeomorphisms $f$ 
which satisfy $f(x+1) = f(x) + 1$
in such a way that its central $\Z$ subgroup corresponds to translations
$T_n: x \mapsto x + n, n \in \Z$. Since $M$ is irreducible,
$H_2(\pi_1(M)) \cong H_2(M)$, while $H_2(M) \cong 0$ since $M$ is a
$\Q$-homology $3$-sphere. Hence the previous lemma implies that
the restriction of $\rho$ to $\pi_1(\hat M)$ lifts to a homomorphism
$\pi_1(\hat M)  \to \widetilde{Homeo}_+(S^1) \subset Homeo_+(\R)$. Since
$\rho$ has non-abelian image, the image of the lifted homomorphism will 
not be the trivial group.
Theorem \ref{homtoLO} now implies the desired conclusion.
\qed

The following corollary will not be used later in the paper,
and we refer the interested reader to \cite{CD} for a definition
of \emph{taut} foliations.

\begin{cor} \label{cal}
{\rm (Calegari-Dunfield)}
Let $M$ be an irreducible, atoroidal $\Q$-homology $3$-sphere
which admits a transversely orientable taut foliation. If $\hat M$ is 
the cover of $M$ corresponding to the abelianization of $\pi_1(M)$, 
then $\pi_1(\hat M)$ is left-orderable.
\qed
\end{cor}

\pf Calegari and Dunfield prove \cite{CD} that under the 
conditions of the corollary, the fundamental group of $M$ admits a 
faithful representation to $Homeo_+(S^1)$ arising from Thurston's
universal circle construction. The corollary therefore follows from 
proposition \ref{hom3}. 
\qed

\begin{cor}
Let $M$ be a Seifert fibred manifold which is also a $\Z$-homology
$3$-sphere. Then $M$ is either homeomorphic to the Poincar\'e homology sphere
$($with $\pi_1(M)$ finite and non-trivial$)$, or else $\pi_1(M)$ is
left-orderable.
\end{cor}

\pf If $M$ is a Seifert fibred manifold and a $\mathbb Z$-homology
$3$-sphere other than the Poincar\'e homology-sphere, it is
$P^2$-irreducible and either it is homeomorphic to $S^3$ or its base
orbifold is hyperbolic (cf. \S 4). The bi-orderability of $\pi_1(M)$ is
obvious in the former case, while in the latter
we observe that the quotient of $\pi_1(M)$ by its centre is a non-trivial
Fuchsian subgroup of $PSL_2(\R) \subset Homeo_+(S^1)$. In this case apply the
previous proposition.
\qed

\begin{exam}
{\rm We shall illustrate proposition \ref{hom3} with the following example. 
Let $M_K$ denote the exterior of the figure $8$ knot $K$. For each
extended rational number $\frac{p}{q} \in \mathbb Q \cup \{\frac10\}$ 
let $M_K(\frac{p}{q})$ be the $\frac{p}{q}$-Dehn filling of $M_K$, that 
is $M_K(\frac{p}{q})$ is the manifold obtained by attaching a solid torus 
$V$ to $M_K$ in such a way that the meridian
of $V$ wraps $p$ times meridionally around $K$ and $q$ times longitudinally.
Each of these manifolds is irreducible and is a $\mathbb Q$-homology
$3$-sphere if and only if $\frac{p}{q} \ne 0$.  We will show that for 
$-4 < \frac{p}{q} < 4$, $\pi_1(M_K(\frac{p}{q}))$ admits a
representation to $PSL_2(\mathbb R)$ with non-abelian image and hence is 
virtually left-orderable (when $\frac{p}{q} = 0$ apply 
Corollary \ref{gencond}). We remark that Dunfield and Thurston \cite{DT} 
have proven that each $M_K(\frac{p}{q})$,
$\frac{p}{q} \ne \infty$, has a finite cover with a positive first Betti number and so it follows that each Dehn filling
of $M_K$ has a virtually left-orderable fundamental group.

There is a presentation of the form
$$\pi_1(M_K) = \langle x,y \; | \; wx = yw \rangle$$
where $x$ represents a meridian of $K$ and $w= xy^{-1}x^{-1}y$.
Given $s \geqslant\frac{1+ \sqrt{5}}{2}$ set 
$t = \frac{1}{2(s-s^{-1})} (1 + \sqrt{(s-s^{-1})^4 + 2(s-s^{-1})^2-3}) 
\in \mathbb R$. The reader can verify that there is a representation
$\phi_s: \pi_1(M_K) \to SL_2(\mathbb R)$ such that
$$\phi_s(x) =  \left(\begin{array}{cc} s & 0 \\ 0 & s^{-1} \end{array}
\right), \;\;\;
\phi_s(y) =
\left(\begin{array}{cc} (\frac{s+s^{-1}}{2}) + t &(\frac{s-s^{-1}}{2}) + t \\
(\frac{s-s^{-1}}{2}) - t &  (\frac{s+s^{-1}}{2}) - t
\end{array} \right).$$
It is simple to see that each $\phi_s$ has a non-abelian image in
$PSL_2(\mathbb C)$
and that $\phi_s$ is reducible if and only if $s = \frac{1 + \sqrt{5}}{2}$.

The elements $\mu = x$ and $\lambda = yx^{-1}y^{-1}x^2y^{-1}x^{-1}y$ of $\pi_1(M_K)$ represent meridian and longitude classes of the knot $K$.
Set $A_s =\phi_s(\mu),  B_s =\phi_s(\lambda)$.
As $A_s$ is diagonal but not $\pm I$ and $[A_s, B_s] = I$, $B_s$
is also diagonal. Let $\zeta(A_s), \zeta(B_s)
\in \mathbb R$ be the $(1,1)$-entries of $A_s, B_s$. Then $\zeta(A_s) = s$
while direct calculation yields
$$\zeta(B_s) = \frac{1}{2s^4} ((s^8 -s^6 -2s^4 -s^2 + 1) + (s^4 -
1)\sqrt{s^8 -2s^6 -s^4 -2s^2 + 1})$$
Now $\phi_s$ induces a homomorphism $\pi_1(M_K(p/q))
\to PSL_2(\mathbb R)$ if and only if $\zeta(A_s)^p\zeta(B_s)^q = \pm 1$, or
equivalently
$$-\frac{ln|\zeta(B_s)|}{ln|\zeta(A_s)|} = \frac{p}{q}.$$
Thus we must examine the range of the function
$$g:[\frac{1 + \sqrt{5}}{2}, \infty) \to \mathbb R, \;\;\; s \mapsto
-\frac{ln|\zeta(B_s)|}{ln|\zeta(A_s)|}.$$
Since $\phi_{\frac{1+ \sqrt{5}}{2}}$ is reducible and $\lambda$ lies in the
commutator subgroup of
$\pi_1(M_K)$, $\zeta(B_{\frac{1+ \sqrt{5}}{2}}) = 1$ and therefore
$ln|\zeta(B_{\frac{1+ \sqrt{5}}{2}})| = 0$.
On the other hand $\zeta(A_{\frac{1+ \sqrt{5}}{2}}) = \frac{1+ \sqrt{5}}{2}
>  1$ so that
$ln|\zeta(A_{\frac{1+
\sqrt{5}}{2}})| > 0$. It follows that $g(\frac{1+ \sqrt{5}}{2}) = 0$.

Next observe that
$\lim_{s \to \infty} \zeta(A_s)^{-4}\zeta(B_s) = \lim_{s \to \infty}
s^{-4}\zeta(B_s) = 1$. Therefore
$$\lim_{s \to \infty} (-4 ln|\zeta(A_s)| + ln|\zeta(B_s)|) = 0,$$
which yields $\lim_{s \to \infty} g(s) = 4$.
Hence the range of $g$ contains
$[0, 4)$ and so for each rational $\frac{p}{q}$ in this interval, there is
at least one $s(\frac{p}{q}) \in (\frac{1 + \sqrt{5}}{2}, \infty)$ such
that $\phi_{s(\frac{p}{q})}$
factors through $\pi_1(M_K(\frac{p}{q}))$. Further the image of this
representation is
non-abelian. Our argument is completed by observing that the amphicheirality of
$M_K$ implies that if $\pi_1(M_K(\frac{p}{q}))$ admits a non-abelian
representation to
$PSL_2(\mathbb R)$, then so does $\pi_1(M_K(-\frac{p}{q}))$.}
\end{exam}


\section{Seifert fibre spaces} \label{sfs}

In this section we develop some background material on Seifert fibred
spaces which will be used later in the paper.  This important class of
3-manifolds was introduced by Seifert \cite{Seif} in 1933, and later extended
to include singular fibres which reverse orientation.  We adopt
the more general definition, as in Scott \cite{Sc2}.
A Seifert fibred space is a $3$-manifold $M$ which is foliated by circles.
It is assumed that each leaf $C$, called a {\it fibre}, has a closed tubular
neighbourhood $N(C)$ consisting of fibres. If $C$ reverses orientation in $M$,
then $N(C)$ is a {\it fibred solid Klein bottle}. A specific model
is given by $$(D^2 \times I)/\{(x, 1) = (r(x), 0)\}$$ where $D^2 \subset
{\mathbb C}$ is the unit disk, $r\co D^2 \to D^2$ is a reflection (e.g.\
complex conjugation), and the foliation is induced from the $I$-factors in
$D^2 \times I$. Note that most fibres wind twice around $N(C)$, but
there is also an annulus consisting of {\it exceptional} fibres,
each of which winds around $N(C)$ once.

If $C$ preserves orientation, then $N(C) \cong S^1 \times D^2$ is a {\it
fibred solid torus}. In this
case the fibre preserving homeomorphism classes of such objects are
parameterised by an integer $\alpha \geqslant
1$ and the $\pm$ class (mod $\alpha$) of an integer $q$ coprime with
$\alpha$. Specific models are represented by
$$(D^2 \times I) / \{(x,1) = (e^{\frac{2\pi i q}{\alpha}}x, 0)\}$$
endowed with the foliation by circles induced from the $I$-factors. The
fibre $C_0$ corresponding to $\{0\}
\times I$ winds once around $N(C)$, while the others wind $\alpha$ times.
If $C = C_0$ we define the {\it index} of $C$ to be $\alpha$, otherwise $1$.
Note that the index of an orientation preserving fibre $C$ is well-defined.
Such a fibre is referred to
as {\it exceptional} if its index is larger than $1$.

The reader will verify that the
space of leaves in $N(C)$ is always a $2$-disk, and therefore the space of
leaves in
$M$, called the {\it base space}, is a surface $B$. There is more structure
inherent in $B$, however.
Indeed, it is the underlying space of a $2$-dimensional orbifold ${\cal
B}$, called the {\it base orbifold}
of $M$, whose singular points correspond to the exceptional
fibres $C$ of the given Seifert structure.
If $C$ preserves orientation, then the associated point in $B$ is
a cone point, lying in $\mbox{int}(B)$, whose order equals the index of
$C$.
If $C$ reverses orientation, then it corresponds to a reflector point
in $\partial B$, which in turn lies on a whole curve of reflector points in
$B$.   The base space $B$ will also be written $|{\cal
B}|$.
There is a short exact sequence (see, for instance, lemma 3.2 of \cite{Sc2})
$$1 \to K \to \pi_1(M) \to \pi_1^{orb}({\cal B}) \to 1 \eqno {(4.1)}$$
where $K$ is the cyclic subgroup of $\pi_1(M)$ generated by a regular fibre
and
$\pi_1^{orb}({\cal B})$ is the orbifold fundamental group of ${\cal B}$
(\cite{Th1}, Chapter 13).

In the case that $M$ is orientable, the singularities of ${\cal B}$ are cone
points lying in the
interior of $B$. We shall say that ${\cal B}$ is of the form $B(\alpha_1,
\alpha_2, \ldots , \alpha_n)$
where $\alpha_1, \alpha_2, \ldots , \alpha_n \geqslant 2$ are the indices
of the
exceptional fibres. Note that
in this case $\partial M$ is foliated by regular fibres and so consists of
tori.

Following are some well-known facts about Seifert fibred spaces which will
be useful.

\begin{prop} \label{sfsbasic}
Suppose $M$ is a compact, connected Seifert fibred space and denote by
$h \in \pi_1(M)$ a class corresponding to a regular $($i.e.
non-exceptional$)$ fibre. \\
$(1)$ If $h$ has finite order, then $M$ is orientable
and finitely covered by $S^3$. In particular, $\pi_1(M)$ is finite. \\
$(2)$ If $h = 1$, then $M \cong S^3$. \\
$(3)$ If $M$ is reducible, then $M = S^1 \times S^2$ or $S^1
\tilde \times S^2$ or $P^3 \# P^3$.  The first two have
$($bi-orderable$)$ group $\Z$ and first
Betti number $b_1(M) = 1$.
However $P^3 \# P^3$ has group $\Z/2 * \Z/2$ $($which is not
left-orderable$)$, first Betti number $0$, and base orbifold $P^2$.\\
$(4)$ If $M$ is nonorientable with 2-torsion in $\pi_1(M)$,
then $M = P^2 \times S^1$ with base orbifold $P^2$.
Its group $\Z/2 \times \Z$ is not left-orderable. \\
$(5)$ If  $\pi_1(M) \cong \Z$, then $M = S^1 \times S^2$ or $S^1
\tilde \times S^2$ or a solid torus or solid Klein bottle.\\
\end{prop}

\pf (1) Let $p:\tilde M \to M$ be the universal cover of $M$. If $h$ has finite
order in $\pi_1(M)$, then the inverse image of each fibre in $M$ is
a circle in $\tilde M$. In this way there is a Seifert
fibring of $\tilde M$ with base orbifold $\tilde {\cal B}$ say and a
commutative diagram
$$ \begin{array}{ccc}
\tilde M & \longrightarrow & \tilde {\cal B} \\
\downarrow &               & \downarrow \\
M &        \longrightarrow & {\cal B}
\end{array} $$
where $\tilde {\cal B} \to {\cal B}$ is an orbifold covering.
The simple-connectivity of $\tilde M$ implies that
$\pi_1^{orb}(\tilde {\cal B}) =\{1\}$ (cf. exact sequence (4.1))
and therefore Riemann's uniformization theorem and theorem 2.3 of
\cite{Sc2} imply that $\tilde {\cal B}$ is
either a contractible surface without cone points  or one of $S^2, S^2(p)$
or $S^2(p, q)$
where $\mbox{gcd}(p,q) = 1$.
The first case is ruled out
as otherwise $\tilde M \cong |\tilde {\cal B}| \times S^1 \simeq S^1$ is not
$1$-connected. In the latter three cases,
$\tilde M$ is a union of two solid tori and therefore must be the
$3$-sphere. Hence the fundamental group of $M$ is finite. Since $S^3$
admits no
fixed-point free orientation reversing homeomorphism, $M$ is orientable.

(2) Next assume that $h = 1$. By (1)  $\tilde M \cong S^3$ and
$\tilde {\cal B}$ is either $S^2, S^2(p)$ or $S^2(p,q)$ where
$\mbox{gcd}(p,q) = 1$.
We also know that $\pi_1(M)$ is finite and $M$ is orientable.
By hypothesis, the inclusion of each fibre of $M$ lifts to an inclusion of
the fibre in
$\tilde M$. It follows that $\pi_1(M)$ acts freely on the components of the
inverse image of any
fibre of $M$. Thus the induced action of $\pi_1(M)$ on $|\tilde {\cal B}|
\cong S^2$ is free and
therefore $\pi_1(M)$ is a subgroup of $\Z/2$. We will show that
$\pi_1(M) \not \cong \Z/2$.

Assume otherwise and observe that since $\pi_1(M)$ freely permutes the
exceptional
fibres of
the Seifert structure on $\tilde M$, the only possibility is for $\tilde
{\cal B} = S^2$.
Exact sequence (4.1) yields $\pi_1^{orb}({\cal B}) \cong \pi_1(M) \cong \Z/2$
and so
$M$ is a  locally trivial $S^1$-bundle over $P^2$.
Splitting ${\cal B}$ into the union of a M\"{o}bius band and a
$2$-disk shows that $M$ is a Dehn filling of the twisted $I$-bundle over the
Klein bottle. A
homological calculation then shows that the order of $H_1(M)$ is divisible
by $4$. But this
contradicts our assumption that $\pi_1(M) \cong \Z/2$. Thus $\pi_1(M) =
\{1\}$ and so $M = \tilde M /\pi_1(M) = \tilde M \cong S^3$.

(3) Suppose that $M$ is reducible and let $S \subset M$ be an essential
$2$-sphere.
The universal cover $\tilde M$ of $M$ is also reducible as otherwise
a $3$-ball bounded by an innermost lift of $S$ to $\tilde M$ projects to a
ball bounded by $S$.
Now the interior of the universal cover of a Seifert fibred space is either
$S^3, \R^3$ or
$S^2 \times \R$ (see eg. [Sc2, Lemma 3.1]) and therefore the interior of
$\tilde M$ is homeomorphic
to $S^2 \times \mathbb R$. By part (1), $h$ has infinite order in $\pi_1(M)$,
and it is not hard to
see that the quotient of $S^2 \times \mathbb R$ by the action of some power 
of $h$ is $S^2 \times S^1$. Thus $M$
itself is finitely covered by $S^2 \times S^1$ and so is one of $S^1 \times
S^2, P^3 \# P^3, S^1
\tilde{\times} S^2$, or $S^1 \times P^2$ [Tlf]. If $M \cong S^1 \times
P^2$, then the fact that
$H_2(S^1 \times P^2) = 0$ implies that $S$ is separating and consideration of
$\pi_1(S^1 \times P^2) \cong \mathbb Z \oplus \mathbb Z/2$ implies that it
bounds a
simply-connected submanifold $A$ of $S^1 \times P^2$. Hence $A$ lifts to
$\tilde M \cong S^2
\times \mathbb R \subset \mathbb S^3$. It follows that
$A$ is a $3$-ball and therefore $M \not \cong S^1 \times P^2$.

(4) Suppose that $M$ is nonorientable with 2-torsion in $\pi_1(M)$ and let
$\tilde M$ be its
universal cover. The group $\pi_1(M)$ is infinite by part (1) and so
$\tilde M$ is non-compact.
In particular $H_q(\tilde M) = 0$ if $q \geqslant 3$.
If $\pi_2(M) = 0$, then $H_q(\tilde M) = 0$ for all
$q$ and so $\tilde M$, being simply-connected, is contractible.  But then the
quotient of $\tilde M$ by a cyclic group of order two $\Z/2 \subset \pi_1(M)$
would be a $K(\Z/2,1)$, which is impossible as $\Z/2$ has infinite
cohomological dimension. Hence $\pi_2(\tilde M) = \pi_2(M) \ne 0$, which
implies that
$\tilde M \cong S^2 \times \mathbb R$ and $M$ is closed (cf. the proof of
part (3)).
Amongst the four closed manifolds covered by $S^2 \times \R$ only  $P^2
\times S^1$ satisfies the
hypotheses of (4).

(5) Suppose that $\pi_1(M) \cong \mathbb Z$. If $\partial M \ne \emptyset$
it contains a
compressible torus or Klein bottle. By parts (3) and (4) $M$ is $\mathbb
P^2$-irreducible and
therefore $M$ is either a solid torus or a solid Klein bottle. On the other
hand if
$\partial M = \emptyset$ and $M$ is orientable, any non-separating closed,
connected, orientable surface
in $M$  (which exists since $b_1(M) = 1$) may be
compressed down to a non-separating $2$-sphere. Thus by part (3) $M$ is
$S^1 \times S^2$.
This implies that if $\partial M = \emptyset$ and $M$ is non-orientable, 
then the orientation double cover of $M$ is $S^1 \times S^2$, so that we 
have $M \cong S^1 \tilde{\times} S^2$ (cf. the argument in part (3)).
\qed

Consider a closed, connected, oriented Seifert manifold $M$.
A useful notation for such manifolds appears in \cite{EHN} which we
describe next. 

The base orbifold of $M$ is of
the form $B(\alpha_1, \alpha_2, \ldots , \alpha_n)$ where $B$ is a closed
surface and $\alpha_1, \alpha_2,
\ldots , \alpha_n \geqslant 2$. As is well-known, $B$ is determined by
$$g = \left\{ \begin{array}{ll} 1 - \frac{\chi(B)}{2} & \mbox{if $B$ is
orientable} \\
\chi(B) - 2 & \mbox{if $B$ is non-orientable.} \end{array} \right. $$
When $n = 0$, $p\co M \to B$ is an
$S^1$-bundle whose total space is oriented, and so $M$ is completely
determined by $g$ and an integer $b$, essentially the Euler number of 
the circle bundle $M \to B$,
measuring the obstruction to the existence of a cross-section (see 
the discussion on pages 434--435 of \cite{Sc2}). An explicit
description of $b$ is obtained as follows. Let $D \subset \mbox{int}(B)$ 
be a $2$-disk and set $B_0 = B \setminus \mbox{int}(D), M_0 = p^{-1}(B_0)$, 
so that $M$ is constructed from $M_0$ by attaching the solid
torus $p^{-1}(D)$, i.e. $M$ is obtained from $M_0$ by a Dehn filling.
The bundle $p_0\co M_0 = p^{-1}(B_0) \to B_0$ is uniquely determined by
the fact that its total space is orientable,  and it can be shown that
$p_0$ admits a section $s$. The orientation on
$M$ determines  orientations on $s(\partial B_0)$ and a circle fibre $H$ on
$p_0^{-1}(\partial B_0)$, and hence a homology basis $[s(\partial B_0)],
[H]$ for $H_1(p_0^{-1}(\partial B_0))$, well-defined up to a simultaneous
change of sign.
Then $b$ is the unique integer  such that the meridional slope of
$p^{-1}(D)$ corresponds to
$\pm ([s(\partial B_0)]  + b [H]) \in H_1(p_0^{-1}(\partial B_0))$.

When $n > 0$ we proceed similarly.
Let $C_1, C_2, \ldots , C_n$ be the exceptional fibres in $M$, $C_0$ a
regular fibre, and $x_0, x_1, x_2,
\ldots , x_n \in B$ the points to which they correspond.  Choose disjoint
$2$-disks $D_0, D_1,D_2, \ldots , D_n \subset B$  where
$x_j \in \mbox{int}(D_j)$ and set $B_0= B \setminus \bigcup_i \mbox{
int}(D_i), M_0 = p^{-1}(B_0)$. The
$S^1$-bundle $M_0 \to B_0$ admits a section $s$ and as in the last
paragraph, the homology classes
of the meridional slopes of the solid tori $p^{-1}(D_j)$ are of the form
$\pm(\alpha_j [s(\partial D_j)] + \beta_j [H_j])$ where
$\alpha_j$ is the index of $C_j$. In fact there is a unique choice
of $s$, up to vertical homotopy, satisfying $0 < \beta_j < \alpha_j$ for
$j = 1, 2, \ldots , n$. Make this choice and set $b = \beta_0$. Then $M$
both determines and is determined
by the integers $g, b$ and the rational numbers $\frac{\beta_1}{\alpha_1},
\frac{\beta_2}{\alpha_2}, \ldots ,
\frac{\beta_n}{\alpha_n} \in (0,1)$. Conversely given such a sequence of
numbers we may construct a closed,
connected, oriented Seifert manifold which realizes them. In the notation
of \cite{EHN},
$$M = M(g; b, \frac{\beta_1}{\alpha_1}, \frac{\beta_2}{\alpha_2}, \ldots ,
\frac{\beta_n}{\alpha_n}).$$
The fundamental group of this manifold is given by
$$\begin{array}{ll}\pi_1(M)= \langle a_1,
b_1,...,a_g,b_g,\gamma_1, \ldots , \gamma_n,h\;|\;h\;\;central,\;\;
\gamma_j^{\a_j}=h^{-\b_j} \;\;
(j=1,...,n)\\
\;\;\;\;\;\;\;\;\;\;\;\;\;\;\;\;\;\;\;\;\;\;\;\;\;\;\;\;\;\;\;\;\;
[a_1,b_1]...[a_g,b_g]\gamma_1...\gamma_n = h^b \rangle
\;\;\;\;\; \end{array}$$
when $g \geqslant 0$, and
$$\begin{array}{ll}\pi_1(M)= \langle a_1,...,a_{|g|},\gamma_1, \ldots ,
\gamma_n,h\;|\;a_jha_j^{-1}=h^{-1}
\;\;\;  (j=1,...,|g|), \;\;\; \gamma_j^{\a_j}=h^{-\b_j}, \\
\;\;\;\;\;\;\;\;\;\;\;\;\;\;\;\;\;\;\;\;\;\;\;\;\;\;\;\;\; \gamma_j h
\gamma_j^{-1}=h \;\; (j=1,...,n),
\;\; a_1^2...a_{|g|}^2\gamma_1... \gamma_n = h^b \rangle.
\end{array}$$
when $g < 0$ (\cite{Jc}, VI.9-VI.10).
The element $h \in \pi_1(M)$ which occurs in these presentations
is represented by any regular fibre of the Seifert
structure. It generates
a normal cyclic subgroup $K$ of $\pi_1(M)$ which is central
if $B$ is orientable.

Let $\chi (B)$ be the Euler characteristic of $B$
and recall that the orbifold Euler characteristic (\cite{Th1}, Chapter 13)
of the orbifold ${\cal B}$ is the rational number given by
$$\chi^{orb}({\cal B})=\chi (B)-\Sigma_{i=1}^{n}(1-\frac{1}{\a_i})
= \left
\{\begin{array}{l}2-2g-\Sigma_{i=1}^{n}(1-\frac{1}{\a_i})\;\;\;\mbox{if $B$
is orientable}\\
2-g-\Sigma_{i=1}^{n}(1-\frac{1}{\a_i})\;\;\;\mbox{if $B$ is non-orientable}.
\end{array} \right.$$
The orbifold ${\cal B}$ is called {\it hyperbolic}, respectively {\it
Euclidean},
if it admits a hyperbolic, respectively Euclidean, structure and this
condition is shown to be equivalent to the condition $\chi^{orb}({\cal B})<
0$, respectively $\chi^{orb}({\cal B})= 0$, in \cite{Th1}, chapter 
13. As such structures are developable, it follows that 
$\pi_1^{orb}({\cal B})$ acts properly discontinuously on $\mathbb 
E^2$ (when $\chi^{orb}({\cal B})= 0$) and on $\mathbb H^2$ (when 
$\chi^{orb}({\cal B}) < 0$) with quotient ${\cal B}$ (\cite{Th1}, 
chapter 13).


\section{Left-orderability and foliations} \label{cod1}

In this paragraph we shall focus on a different class of objects,
namely on codimension 1 foliations. Such foliations 
will play an important role in our analysis of the left-orderability 
of the fundamental groups of Seifert fibred manifolds. Their 
connection to orderabilty comes through the induced action of the 
fundamental group of the ambient manifold on the leaf space of the 
induced foliation on the universal cover. Under certain natural hypotheses 
this leaf space can be shown to be homeomorphic to the real line.  
Throughout, the foliations we will 
consider will be ${\cal C}^1$ and transverse to the boundary of the 
ambient manifold. 

We saw in theorem \ref{eff} that a countable 
group $G$ is left-orderable
if and only if it acts effectively on $\mathbb R$ by order-preserving
homeomorphisms. In the case of the fundamental group of a $P^2$-irreducible
$3$-manifold $M$, theorem \ref{homtoLO} shows that this condition can be
relaxed to the existence of a homomorphism $\pi_1(M) \to 
Homeo_+(\mathbb R)$ with
non-trivial image. Given such a homomorphism, our next lemma shows 
that it can be supposed to induce a {\it non-trivial action} on 
$\mathbb R$, that is, the action  has no global fixed point.

\begin{lemma} \label{obs}
If there is a homomorphism $G \to Homeo_+(\mathbb R)$ with image $\ne \{id\}$,
then there is another such homomorphism which induces an action on 
$\mathbb R$ without
global fixed points.
\end{lemma}

\pf Fix a homomorphism $\phi\co G \to Homeo_+(\mathbb R)$ with image
$\ne \{id\}$ and observe that
$F := \{x \; | \; \phi(\gamma)(x) = x \mbox{ for every }\gamma \in G\}$ is 
a closed, proper subset of $\mathbb R$. Each component $C$ of the non-empty 
set $\mathbb R \setminus F$
is homeomorphic to $\mathbb R$ and is invariant under the given action.
By restricting the action to $C$, we obtain the desired
action without global fixed points.  \qed

We can therefore suppose, when necessary, that if we have a 
homomorphism of 
a group $G$ to $Homeo_+(\mathbb R)$ which has a 
non-trivial image, the associated action on $\mathbb R$ is 
non-trivial. 

Gabai raised the problem of developing a theory of 
non-trivial group actions on order trees and asked some fundamental 
questions about the nature of those $3$-manifolds whose groups admit 
such actions, especially those which act on $\mathbb R$ (\S 4, 
\cite{Ga}). Using standard techniques, it is possible to translate 
the existence of such  actions into a topological condition. Indeed, 
if $M$ is a compact, connected, orientable $3$-manifold, we have just 
seen that a necessary and sufficient condition for $\pi_1(M)$ to be 
left-orderable is that there be a non-trivial action of $M$ 
determined by some homomorphism $\phi\co \pi_1(M) \to Homeo_+(\mathbb 
R)$. Given such a homomorphism, one can construct (cf. remark 4.2 
(ii), \cite{Ga}) a transversely orientable, transversely essential 
lamination whose order tree maps $\pi_1(M)$-equivariantly, with 
respect to $\phi$, to $\mathbb R$. As it will not play any subsequent 
role in the paper, we direct the reader to \cite{Ga} for definitions 
and details.

One way to produce actions of a $3$-manifold group 
$\pi_1(M)$ on the reals is by constructing $\mathbb R$-covered foliations. 
These are codimension $1$ foliations such that the space of leaves
of the pull-back foliation in the universal cover of $M$ is $\mathbb R$.
Many examples of
hyperbolic $3$-manifolds with $\mathbb R$-covered foliations exist. See
\cite{Fe}, \cite{Th2}, and \cite{Ca1,Ca2} for
various constructions and related information. See, however,
section \ref{hyp} for examples of hyperbolic 3-manifolds which do not 
contain such foliations.

Recall that a codimension one foliation of a 3-manifold
is {\it orientable} if its 
tangent field is orientable. It is {\it transversely orientable} if 
its transverse line field is orientable. Clearly these notions are 
equivalent if and only if the ambient manifold is orientable. 

\begin{lemma}\label{actleafspace}
Let $M$ be a compact, connected $3$-manifold and
${\cal F}$ a transversely oriented, $\mathbb R$-covered foliation in 
$M$. Denote by $\tilde {\cal F}$ the lift of ${\cal F}$ to
$\tilde M$ and let $\phi\co \pi_1(M) \to Homeo(\mathbb R)$
be the homomorphism induced by the action
of $\pi_1(M)$ on $\tilde {\cal F}$. Then the image of $\phi$ lies in 
$Homeo_+(\mathbb R)$. 
\end{lemma}

\pf The transverse field to $\tilde {\cal F}$ is the pull-back of 
that of ${\cal F}$ and so is orientable. It is easy to see that the 
group of deck transformations preserves either of its orientations. 
Thus the lemma holds. 
\qed

\begin{prop} \label{fol}
Let $M$ be a compact, connected, $P^2$-irreducible $3$-manifold which
admits a transversely oriented, $\mathbb R$-covered  foliation ${\cal 
F}$ of $M$. Then the fundamental group of $M$ is left-orderable.
\end{prop}

\pf First note that there is a compact subset $C$ of the leaf space 
$\mathbb R$ of $\widetilde{{\cal F}}$ which meets each orbit of the action 
of $\pi_1(M)$ on $\mathbb R$. (Take $C$ to be the image in $\mathbb R$ 
of any fundamental domain of the universal cover $\tilde M \to M$.) 
Thus the homomorphism $\pi_1(M) \to Homeo(\mathbb R)$ associated to 
the action has a non-trivial image. Lemma \ref{actleafspace} shows 
that its image lies in $Homeo_+(\mathbb R)$ and thus  theorem 
\ref{homtoLO} implies that $\pi_1(M)$ is LO. 
\qed

For Seifert manifolds, there is a distinguished class of 
codimension $1$ foliations. \\

\noindent {\bf Definition:} A {\it horizontal foliation} of a Seifert
fibred manifold is a foliation of $M$ by (possibly noncompact) 
surfaces which are everywhere transverse to the Seifert fibres. \\ 

Though such foliations are traditionally referred to as {\it 
transverse}, we have chosen to use the equally appropriate term {\it 
horizontal} to avoid confusion with the notion of a transversely 
oriented foliation discussed in the next section. It is shown in 
corollary 4.3 of \cite{EHN} that a Seifert fibred manifold which 
admits a ${\cal C}^0$ horizontal foliation also admits an analytic 
horizontal foliation. 

The combined work of various authors has 
resulted in a complete understanding of which Seifert bundles admit 
horizontal foliations. In the following theorem we consider the case 
where $M$ is closed and $g = 0$.

\begin{thm} \label{class} {\rm (\cite{EHN}, \cite{JN2}, \cite{Na})}
Let $M = M(0;b, \frac{\beta_1}{\alpha_1}, \ldots ,
\frac{\beta_n}{\alpha_n})$ be an orientable Seifert
fibred manifold   where $n \geqslant 3$, $b \in \mathbb Z$ and $\alpha_j,
\beta_j$ are integers for which
$0 < \beta_j < \alpha_j$. Then $M$ admits a horizontal foliation if and
only if one of the following conditions holds: \\
$(1)$ $-(n-2) \leqslant b \leqslant -2$. \\
$(2)$ $b = -1$ and there are relatively prime integers $0 < a < m$ such
that for some permutation
$(\frac{a_1}{m}, \ldots , \frac{a_n}{m})$ of 
$(\frac{a}{m}, \frac{m-a}{m},\frac{1}{m},\ldots ,
\frac{1}{m})$ we have $\frac{\beta_j}{\alpha_j} < \frac{a_j}{m}$ for each
$j$.  \\
$(3)$ $b = -(n-1)$ and after replacing each $\frac{\beta_j}{\alpha_j}$ by
$\frac{\alpha_j - \beta_j}{\alpha_j}$, condition $(2)$ holds.
\qed
\end{thm}

Roberts and Stein have shown \cite{RS} that a necessary 
and sufficient condition for the fundamental group of an irreducible, 
non-Haken Seifert fibred manifold to act non-trivially on $\mathbb R$ 
is that the manifold admit a horizontal foliation dual to the action.
We shall offer a new proof, and expand on this theme in this section and 
the next.  

A horizontal foliation in a Seifert fibred manifold is orientable if 
and only if the base orbifold is orientable. This follows from the 
observation that away from exceptional fibres, the tangent field to 
the foliation is the pull-back of the tangent bundle of the orbifold. 
A horizontal foliation in a Seifert fibred manifold is transversely 
orientable if and only if the circle fibres may be coherently 
oriented, or equivalently, there are no vertical Klein bottles in the 
ambient $3$-manifold. Thus, 

\begin{lemma} \label{transorien}
Let $M$ be a compact, connected, orientable Seifert fibred manifold and
${\cal F}$ a horizontal
foliation in $M$. Then ${\cal F}$ is transversely orientable if and only if
the surface underlying the base orbifold of $M$ is orientable. 
\qed
\end{lemma}

The next lemma is well-known. We include its proof for completeness.

\begin{lemma} \label{slr}
Let $M$ be a closed, connected, $P^2$-irreducible Seifert fibred
manifold with infinite fundamental group. Let ${\cal F}$ be a horizontal
foliation in $M$ and $\tilde {\cal F}$ the lift of
${\cal F}$ to $\tilde M$, the universal cover of $M$. Then there is a
homeomorphism $\tilde M \to \mathbb R^3$ which sends $\tilde {\cal F}$
to the set of planes parallel to the $x$-$y$-plane.
In particular ${\cal F}$ is $\mathbb R$-covered.
\end{lemma}

\noindent Before sketching a proof of this lemma, we derive the following
consequence of lemma \ref{transorien}, lemma \ref{slr}, and 
proposition \ref{fol} (cf. corollary \ref{cal}).

\begin{prop} \label{hfro}
Let $M$ be a closed, connected, irreducible, orientable Seifert fibred
manifold. Suppose the surface underlying the base orbifold of $M$ is
orientable, and $M$ admits a horizontal foliation. If $\pi_1(M)$ is
infinite, it is left-orderable.
\qed
\end{prop}

\noindent {\bf Proof of Lemma \ref{slr}.}
Since $\pi_1(M)$ is infinite and $M$ is irreducible, the universal orbifold
cover of ${\cal B}$ is $\R^2$. Pulling back the Seifert fibration via this 
orbifold covering shows that there is a regular covering space 
$\hat M \to M$ where $\hat M$ is an $S^1$-bundle over 
$\mathbb R^2$.  Hence $\hat M$ can be identified with 
$\mathbb R^2 \times S^1$ in such a way that the Seifert
circles pull back to the $S^1$ factors, and so $\tilde M$ is 
identifiable with $\R^3$ in such a way that  the Seifert circles 
pull back to the field of lines parallel to the $z$-axis. Note as well 
that if $\tau\co \mathbb R^3 \to \mathbb R^3$ is vertical translation 
by $1$, then $\tau$ may be taken to be a deck transformation
of the universal cover $\mathbb R^3 \to M$. In particular $\tilde {\cal F}$
is invariant under $\tau$. Let
$$p\co \R^3 \to \R^2$$
be the projection onto the first two coordinates. We will show first of 
all that the restriction of $p$ to any leaf of $\tilde {\cal F}$ is a 
homeomorphism.

Fix a leaf $L$ of $\tilde {\cal F}$ and consider $p|L$. That $p|L$ is 1-1
follows from a classic result of Novikov: 
a closed loop which is everywhere transverse to a codimension-1
foliation without Reeb components is not null-homotopic.  (The case when ${\cal F}$ is ${\cal 
C}^2$ is handled in \cite{No}. See the discussion in \cite{Ga1}, p. 
611, for the general case.) Note that horizontal foliations can never contain a Reeb component as the boundary of such a component is disjoint from every circle transverse to the foliation. 
If there are points $(x_0, y_0, z_0), (x_0, y_0, z_1) \in L$ where $z_0 >
z_1$, the vertical path between them
concatenated with a path in $L$ may be perturbed to be everywhere
transverse to $\tilde {\cal F}$ (this uses
the fact that $\tilde {\cal F}$ is transversely oriented).  Since all 
loops in $\R^3$ are contractible, Novikov's result shows that this is 
impossible. Thus $p|L$ is injective.

Surjectivity follows from the fact that $\tilde {\cal F}$ is transverse to
the vertical line field and that
it is invariant under $\tau$. Here is a more detailed argument.
Firstly, transversality implies that
$p(L)$ is open in $\mathbb R^2$.  We claim
that $\mathbb R^2 \setminus p(L)$ is open as well.

Suppose $(x_0, y_0) \in \R^2 \setminus p(L)$ and let $Z_0 \subset \R^3$ denote
the vertical line through this point.  For any $z \in [0,1]$,
transversality implies there is an open neighborhood $U_z \subset \R^3$ of
$(x_0,y_0,z)$ with the property that any leaf of $\tilde {\cal F}$ that
intersects $U_z$ will also intersect $Z_0$.  By compactness, a finite number
of such $U_z$ will cover $(x_0, y_0) \times [0,1]$, and one can find
$\epsilon > 0$ so that $N_{\epsilon}(x_0, y_0) \times [0,1]$ has the same
property. Since $\tilde {\cal F}$ and $Z_0$ are both $\tau$-invariant, it
follows that $N_{\epsilon}(x_0, y_0) \subset \R^2 \setminus p(L)$, and we have verified that $\R^2 \setminus p(L)$ is open.  The
connectivity of $\mathbb R^2$ implies that $p|L$ is onto, and we've 
shown that $p|L$ is a homeomorphism of $L$ onto $\mathbb R^2$ for each leaf
$L \in \tilde {\cal F}$ .

It follows that each leaf of $\tilde {\cal F}$ intersects each vertical 
line in $\mathbb R^3$ exactly once and so the leaf space 
${\cal L}(\tilde {\cal F})$ is homeomorphic to $\mathbb R$. Let 
$f\co \mathbb R^3 \to \mathbb R$ be the composition of the map
$\mathbb R^3 \to {\cal L}(\tilde {\cal F})$ with such a homeomorphism 
and observe that the map $p \times f\co \mathbb R^3 \to \mathbb R^3$ 
defines a homeomorphism which sends $\tilde {\cal F}$ to the
set of horizontal planes, which is what we set out to prove.
\qed


\section{Left-orderability and Seifert fibred spaces} \label{proof}

In this section we prove theorem \ref{SeifertLO}.
That is, for the fundamental group of a compact, connected, Seifert fibred
space $M$ to be left-orderable, it is necessary and sufficient that 
one of the following holds: \\

\indent $-$ $M \cong S^3$; or\\
\indent $-$ $b_1(M) > 0$ and $M \not \cong {P^2} \times S^1$; or \\
\indent $-$ $b_1(M) = 0$, $M$ is orientable, $\pi_1(M)$ is infinite, 
the base orbifold of $M$
is of the \\ \indent \hspace {1mm} form $S^2(\alpha_1, \alpha_2,$ 
$\ldots, \alpha_n)$, and $M$ admits
a horizontal foliation. \\
Throughout we take $M$ to be a compact,
connected Seifert fibred space with base orbifold ${\cal B}$. By 
proposition \ref{sfsbasic} we may suppose that $M$ is 
$P^2$-irreducible and has a non-trivial fundamental
group. 

\subsection{Sufficiency}

Since $M$ is $P^2$-irreducible, theorem \ref{homtoLO} shows that 
$\pi_1(M)$ is LO  when $b_1(M) > 0$. When $b_1(M) = 0$, lemma 
\ref{posb1} shows that $M$ is closed and orientable. Then proposition 
\ref{hfro} shows that it is LO when
$\pi_1(M)$ is infinite, the base orbifold of $M$ is of
the form $S^2(\alpha_1, \alpha_2, \ldots, \alpha_n)$, and $M$ admits a
horizontal  foliation.   

\subsection{Neccesity}

Assume that $\pi_1(M)$ is LO.  If $b_1(M) > 0$ then the 
$P^2$-irreducibility of $M$ implies that $M \not \cong {P^2} \times 
S^1$, so we are done. Assume then that $b_1(M) = 0$. By lemma 
\ref{posb1}, $M$
must be closed and orientable. Further, since $\pi_1(|{\cal B}|)$ is 
a quotient of $\pi_1^{orb}({\cal B}) \cong \pi_1(M)/\langle h 
\rangle$ (cf. \S 4), $H_1(|{\cal B}|)$ is finite. Thus 
${\cal B} = 
S^2(\alpha_1, \ldots , \alpha_n)$ or
$P^2(\alpha_1, \ldots , \alpha_n)$. Note as well that
$\pi_1(M)$ is infinite as it is a non-trivial torsion free group. We 
must prove that ${\cal B} = S^2(\alpha_1, \ldots , \alpha_n)$ and $M$
admits a horizontal foliation.

First observe that if ${\cal B} = 
S^2(\alpha_1, \ldots , \alpha_n)$ then $n \geqslant 3$ and if $n = 3$ then 
$(\alpha_1, \alpha_2, \alpha_3)$ is a Euclidean or hyperbolic triple. 
Otherwise $M$ would be $S^1 \times S^2$ or have a finite fundamental 
group (see \cite{Jc}, VI.11 (c)). Thus $\chi^{orb}({\cal B}) \leq 0$ 
so that ${\cal B}$ admits a Euclidean or hyperbolic structure (cf. \S 
4). In particular $\pi_1^{orb}({\cal B})$ acts properly 
discontinuously on $\mathbb E^2$ or $\mathbb H^2$ with quotient 
${\cal B}$. We also note that such $M$ admit a unique Seifert 
structure up to isotopy (see \cite{Jc}, theorem VI.17). 

When 
${\cal B} = P^2(\alpha_1, \ldots , \alpha_n)$ then $n \geqslant 2$ as 
otherwise $M$ would be either $P^3 \# P^3, S^1 \times S^2$ or have a 
finite fundamental group (see \cite{Jc}, VI.11 (c)). Thus 
$\chi^{orb}({\cal B}) \leqslant 0$ and again we see that ${\cal B}$ admits 
a Euclidean or hyperbolic structure and $M$ admits a unique Seifert 
structures up to isotopy. 

Express $M$ in the form $M(g;b, 
\frac{\beta_1}{\alpha_1}, \ldots ,
\frac{\beta_n}{\alpha_n})$ (cf. \S \ref{sfs})
where $b \in \mathbb Z$,  $\alpha_j, \beta_j$ are integers for
which
$0 < \beta_j < \alpha_j$, and
$$g = \left\{\begin{array}{rl} 0 & \mbox{ when } |{\cal B}| = S^2 \\
-1& \mbox{ when } |{\cal B}| = P^2. \end{array} \right. $$
We noted in \S \ref{sfs} that the
fundamental group of $M$ admits a presentation of the form
$$\left\{\begin{array}{cl}
\langle \gamma_1, \ldots , \gamma_n, h\; | \; h \mbox{ central },
\gamma_j^{\alpha_j} = h^{-\beta_j}, \gamma_1\gamma_2 \ldots
\gamma_n = h^b \rangle  & \mbox{ when } |{\cal B}| = S^2 \\
& \\
\langle \gamma_1, \ldots , \gamma_n, y, h\; | \; 
\gamma_j^{\alpha_j} = h^{-\beta_j}, yhy^{-1}= h^{-1},
y^2\gamma_1\gamma_2 \ldots \gamma_n = h^b \rangle & 
\mbox{ when } |{\cal B}| = P^2, \end{array}
\right. $$
where $h \in \pi_1(M)$ is represented by a regular fibre.  Since 
$\{1\} \ne \pi_1(M)$ is LO, there is a non-trivial homomorphism
$$\phi\co \pi_1(M) \to Homeo_+(\mathbb R).$$
By lemma \ref{obs} we may suppose that the associated action of 
$\pi_1(M)$ on $\mathbb R$ is non-trivial.

\begin{lemma} \label{sh(1)} {\rm (compare lemma 2, \cite{RS})}
Suppose that $M = M(g;b, \frac{\beta_1}{\alpha_1}, \ldots ,
\frac{\beta_n}{\alpha_n})$ where
$g \in \{0, -1\}$. For a homomorphism
$\phi\co \pi_1(M) \to Homeo_+(\mathbb R)$ the following statements are
equivalent: \\
$(1)$ The action induced by $\phi$ is non-trivial. \\
$(2)$ $\phi(h)$ is conjugate in $Homeo_+(\mathbb R)$ to translation by $1$.
\end{lemma}

\pf It is clear that condition (2) implies condition (1), so
suppose the action induced by $\phi$ is non-trivial. As any
fixed point free element of $Homeo_+(\mathbb R)$ is conjugate
to translation by $1$, we shall assume that there is some 
$x_0 \in \mathbb R$ such that $\phi(h)(x_0) = x_0$ and proceed by 
contradiction.  Recall the presentation for $\pi_1(M)$ described
above. We have 
$\phi(\gamma_j)^{\alpha_j}(x_0) = \phi(h)^{-\beta_j}(x_0) = x_0$
for each
$j \in \{1,2, \ldots , n\}$. As $\phi(\gamma_j)$ preserves orientation,
this implies that $x_0$ is fixed by $\gamma_j$. In the case where 
$|{\cal B}| = P^2$ we also have 
$\phi(y)^2(x_0) = \phi(y)^2\phi(\gamma_1\gamma_2\ldots \gamma_n)(x_0) =
\phi(h)^b(x_0) = x_0$ and so $\phi(y)(x_0) = x_0$ as well. In either case
$x_0$ is fixed by $\pi_1(M)$, contradicting the fact that the action is 
non-trivial. Thus $\phi(h)$ is fixed-point
free and therefore is conjugate to translation by $1$.
\qed

Now we complete the proof of theorem \ref{SeifertLO}. By lemma 
\ref{sh(1)} there is a homomorphism 
$\phi: \pi_1(M) \to 
Homeo_+(\mathbb R)$ such that $\phi(h)$ is translation by $1$. We 
noted above that $\pi_1^{orb}({\cal B}) =
\pi_1(M)/\langle h\rangle$ acts properly discontinuously on
$X = \mathbb E^2$ or $\mathbb H^2$. The subsequent diagonal action of
$\pi_1(M)$ on $\mathbb R^3 = X \times \mathbb R$ can be seen to be
free and properly discontinuous. It follows that the quotient $N$ of 
$\mathbb R^3$ by this action 
is a $K(\pi_1(M), 1)$. Since $\pi_1(M)$ 
is infinite, the main result of \cite{Sc3} implies that $M$ is 
is homeomorphic to $N$. The
lines $\{x\} \times \mathbb R$ and the planes $X \times \{t\}$ are
invariant under the action of $\pi_1(M)$ on $X \times \mathbb R$ and
induce, respectively, a Seifert structure on $M$ (necessarily 
isotopic to the one we started with) and a horizontal foliation in 
$M$. Further since the image
of $\phi$ lies in $Homeo_+(\mathbb R)$, an orientation of the
vertical line field in $\mathbb
R^3$ descends to a coherent orientation of the circle fibres in $M$. Thus
the induced foliation is transversely orientable and so $|{\cal B}|$ is
orientable (lemma \ref{transorien}).
It follows that ${\cal B} = S^2(\alpha_1, \ldots , \alpha_n)$.
This completes the proof of theorem \ref{SeifertLO}.
\qed

\begin{rem}
{\rm $\;$ \\
(1) It is proved in \cite{EHN} that there is a homomorphism
$\phi\co \pi_1(M) \to Homeo_+(\mathbb R)$ for which $\phi(h)$ is 
translation by $1$ if and only if $M$ admits a transversely oriented 
horizontal foliation. We have already described
how to construct horizontal foliations from such representations and
conversely how to produce such a
representation when given a horizontal foliation, at least when $b_1(M) =
0$. In particular we have reproved the result of \cite{EHN} in the
special case $b_1(M) = 0$.\\
(2) Lemma \ref{sh(1)} does not hold when $|{\cal B}| \ne S^2, P^2$ and this
explains why the
condition that $\pi_1(M)$ be left-orderable does not imply, in
general, that $M$ admits a horizontal foliation.}
\end{rem}


\section{Bi-orderability and surface groups} \label{biosurf}

All surface groups other than $\Z/2 \cong \pi_1(P^2)$ are locally
indicable and hence LO (cf. theorem \ref{OLI}). To see this, it suffices to
observe that the cover corresponding to a given nontrivial finitely generated
subgroup has infinite torsion-free homology.
Our interest here focuses on the bi-orderability of these groups. We prove,\\

\noindent{\bf Theorem \ref{bisur} }
{\it If $S$ is any connected surface other than the
projective plane $P^2$ or Klein bottle $K = P^2 \# P^2$, then $\pi_1(S)$ is
bi-orderable.
}\\

(Another approach to the bi-orderability of surface groups can 
be found in a recent paper of Champetier and Guirardel \cite{CR})

The theorem is already well-known in the case of orientable surfaces: it is
proved in \cite{Baumslag, Long} that their fundamental groups are
residually free (and hence bi-orderable). However, the fundamental group
of the non-orientable surface $S= P^2 \# P^2 \# P^2$ is not residually free;
this is because the image of any homomorphism $\phi$ from  $\pi_1(S) =
\langle a,b,c \ | \ a^2b^2c^2$ $=1 \rangle$ to a free group is cyclic
(see \cite{LySch}, p.51) and therefore sends the commutator subgroup to
$\{1\}$. For another approach see \cite{GagSp}.

In the remainder of this section we will outline a proof of this theorem.
In fact,
our argument fits into a larger picture, in that similar arguments have been
applied to quite diverse situations - see \cite{RolW} (which contains
further details) as well as \cite{G-M} and \cite{KR}.
In what follows we will denote the connected sum of $n$ projective
planes by $n P^2$.

If $S$ is noncompact, or if $\partial S$
is nonempty, then $\pi_1(S)$ is a free group, and therefore bi-orderable.
Thus we are reduced to considering closed surfaces. According to the
standard classification, such surfaces are either a connected sum of tori,
or projective planes in the nonorientable case.
We remarked above that $\pi_1(P^2)$ is not LO.
For $S = 2P^2$, the Klein bottle, we have already seen in \S \ref{S:Intro}
that $\pi_1(S)$ is LO but not O. The key to our analysis
will be the nonorientable surface with Euler characteristic $-1$, namely
$3P^2$.

\begin{prop} \label{3P}
Let $S = 3P^2$ be the connected sum of three projective planes.  Then
$\pi_1(S)$ is bi-orderable.  \end{prop}

Before proving this result, we explain how it implies theorem \ref{bisur}.
Starting with the nonorientable surfaces $(n+2)P^2 = T^2 \# nP^2$, we
note that $S=3P^2 = T^2 \# P^2$ can be pictured as a torus with a small
disk removed, and replaced by
sewing in a M\"obius band.  Consider an $n$-fold cover of the torus by
itself, and modify the covering by replacing one disk downstairs, and
$n$ disks upstairs, by M\"obius bands.  This gives a covering of $S$ by
the connected sum of $T^2$ with $n$ copies of $P^2$.  Thus the
fundamental group of $(n+2)P^2$ embeds in that of $3P^2$, and is
therefore bi-orderable.

For the orientable surfaces $S_g$ of genus $g\geqslant 2$ (the cases $g=0,1$
being easy) the result follows because $S_g$ is the oriented double cover of
$(g+1)P^2$; so $\pi_1(S_g)$ is a subgroup of a bi-orderable group.
This completes the proof of theorem \ref{bisur}, assuming \ref{3P}.\qed

To prove proposition \ref{3P}, our strategy is to define a surjection from
$G = \pi_1(S)$ to $\Z^2$ with a certain kernel $F$,
so that we have a short exact sequence
$$
1 \longrightarrow F  \longrightarrow \ G \longrightarrow \Z^2
\longrightarrow 1.
$$
Moreover, we shall construct a biordering on $F$ so that the conjugation
action of $G$ on $F$ is by order-preserving automorphisms. By lemma
\ref{oext}, this yields a biordering of $G$.

We recall that $S$ is a torus with a disk removed and a M\"obius band glued
in its place. Squashing that M\"obius band induces the desired surjection
$\psi\co G=\pi_1(S) \to \pi_1(T^2)=\Z^2$. More explicitely, $G$ has
presentation
$$
G = \langle a, b, c : aba^{-1}b^{-1} = c^2 \rangle.
$$
(with $a$ and $b$ corresponding to a free generating set of the punctured
torus, and $c$ corresponding to a core curve of the M\"obius band), and
$\psi$ kills the generator $c$.

The kernel $F$, consisting of those elements
with exponent sums in both $a$ and $b$ equal to zero, is an infinitely
generated free group, with one generator for every element of $\Z^2$.
Geometrically, we can interpret $F$ as the fundamental group of a
covering space $\widetilde{S}$ of $S$: starting with the universal cover
$\R^2 \to T^2$, we remove from $\R^2$ a family of small disks
centered at the integral lattice points, and glue in  M\"obius bands in
their place. Thus we obtain a covering space $\widetilde{S}$ of $S$.

There is no canonical free generating system for $F$ - for definiteness
we may take $$
x_{i,j} = a^ib^jcb^{-j}a^{-i}.
$$
So we have $F = \langle x_{i,j}\rangle; (i,j) \in \Z^2$.

Now $G$ acts upon $F$ by
conjugation, which may be described in terms of the generators as follows.

\begin{lemma} \label{conj}
Suppose $g \in G$ has exponent sums $m$ and $n$ in $a$ and $b$,
respectively.  Then there are $w_{i,j} \in F$ such that
$$gx_{i,j}g^{-1} = w_{i,j}x_{i+m,j+n}w_{i,j}^{-1}.$$
\end{lemma}

\pf Just take $w_{i,j} = ga^ib^{-n}a^{-i-m}$.  Check exponent sums
to verify $w_{i,j} \in F$. \qed

For the following, $F_{ab}$ denotes the abelianization of $F$, which
is an infinitely generated free abelian group, with generators, say
$\widetilde{x}_{i,j}$.
Any automorphism $\varphi$ of $F$
induces a unique automorphism $\varphi_{ab}$ of $F_{ab}$.  For example,
in the above lemma, conjugation by $g$ acts under abelianization as the shift
$\widetilde{x}_{i,j} \to \widetilde{x}_{i+m,j+n}$. Proposition \ref{3P}
now follows from the

\begin{lemma}
There is a bi-ordering of the free group
$F = \langle x_{i,j}\rangle; (i,j) \in \Z^2$ which is invariant under
every automorphism $\varphi\co F \to F$ which induces, on $F_{ab}$, a
uniform shift automorphism $\widetilde{x}_{i,j} \to \widetilde{x}_{i+m,j+n}.$
\end{lemma}

\pf We use the Magnus expansion $\mu\co F \to \Z[[X_{i,j}]]$, where
$\Z[[X_{i,j}]]$ is the ring of formal power series in the infinitely many
noncommuting variables $X_{i,j}$, with the restriction that each power series
may involve only finitely many variables.
The Magnus map $\mu$ is given by
$$ \mu(x_{i,j}) = 1 + X_{i,j};
\quad   \mu(x_{i,j}^{-1}) = 1 - X_{i,j} + X_{i,j}^2 - X_{i,j}^3 + \cdots
$$
Clearly the image of $F$ lies in the group $\Gamma$ of units with constant
term unity,
$\Gamma = \{1 + O(1)\}\subseteq \Z[[X_{i,j}]]$, and the image of the
commutator $[F,F]$ lies in $\{1 + O(2)\}$.
As done in \cite{MKS} for the finitely-generated case, one can prove that
$\mu\co F \to \Gamma$ is an {\it embedding} of groups.
Elements of  $\Z[[X_{i,j}]]$ may be
written in standard form, arranged in ascending degree, and within a
degree terms are arranged lexicographically by their subscripts (which in
turn are ordered lexicographically).  Then two series are compared by the
coefficient of the ``first'' term at which they differ (here is where the
finiteness assumption is necessary). It is well-known (see e.g.\ \cite{KR})
that, restricted to $\Gamma$, this ordering is a (multiplicative) bi-ordering.

Finally, we check that the ordering is invariant under the action by
$\varphi$. Since $\varphi(x_{i,j}) = x_{i+m,j+n}\ c_{i,j}$, where $c_{i,j}$ is
in the commutator subgroup $[F,F]$, and since $[F,F]$ maps into
$\{1 + O(2)\}$ under the Magnus embedding, we have for any $u\in F$ that
the lowest nonzero-degree terms of $\mu(\varphi(u))$ coincide precisely with
those of $\mu(u)$, except that all the subscripts are shifted according
to the rule $X_{i,j} \to X_{i+m,j+n}$. This implies that the Magnus-ordering
of $F$ is invariant under $\phi$. \qed


\section{Bi-orderability and Seifert fibred spaces} \label{bioseif}

Our goal is to prove theorem \ref{SeifertBO}: for the fundamental
group of a compact, connected Seifert fibred space $M$ to be bi-orderable, 
it is necessary and sufficient that it be one of $S^3, S^1 \times S^2, 
S^1 \tilde \times S^2$, a solid Klein bottle, or a locally trivial, 
orientable circle bundle over a surface different from $S^2, P^2$ or
the Klein bottle $K$.


\subsection{Sufficiency}
If $M$ is one of  $S^3, S^1 \times S^2, S^1 \tilde \times S^2$, or a solid
Klein bottle, it is clear that
$\pi_1(M)$ is bi-orderable.  If $M$ is an {\it orientable} circle
bundle over a surface
$B \ne S^2, P^2, K$, then $\pi_2(B)$ is trivial and the
homotopy sequence of the bundle yields the exact sequence:
$$
1 \to \pi_1(S^1) \to \pi_1(M) \to \pi_1(B) \to 1.$$ 
(This sequence 
coincides with that of (4.1) in our present context.)
Since $M \to B$ is an orientable $S^1$-bundle, the bi-orderable group
$\pi_1(S^1)$ is central in $\pi_1(M)$. Theorem \ref{bisur} shows that
$\pi_1(B)$
is bi-orderable, and therefore by lemma \ref{oext}, $\pi_1(M)$ is
bi-orderable as well.


\subsection{Necessity}
Throughout this subsection we use ${\cal B}$ to denote the base orbifold of
$M$,
$B$ to denote the surface underlying ${\cal B}$, $\Sigma \subset B$ to
denote the
singular points of ${\cal B}$, and $L = \partial B \cap \Sigma$ to denote
the set
of reflector
lines of ${\cal B}$.  We are assuming the following:

\bigskip\noindent
($*$)  $M$ is a compact Seifert fibred 3-manifold whose fundamental group
is bi-orderable.
\bigskip

\noindent and our aim is to conclude that $M$ belongs to the given list.

\begin{lemma} \label{orientable} If $M$ satisfies $(*)$, then  \\
$(1)$  the restriction of $M \to B$ to $B \setminus
\Sigma$ is an orientable circle bundle. \\ 
$(2)$ an element of 
$\pi_1(M)$ represented by an arbitrary fibre of $M$ is central.
\end{lemma}

\pf (1) If the bundle in question were not orientable, there would be a 
simple closed
curve $C$ in $B\setminus \Sigma$ over which fibres could not be coherently
oriented.  Then $M$ contains a Klein bottle over $C$. In particular 
$M \not \cong S^3$ and so
by proposition \ref{sfsbasic} (2), the class $h \in \pi_1(M)$ of a 
regular fibre is
non-trivial. If $\gamma \in \pi_1(M)$ corresponds to $C$, then
$1 \ne h^{-1} = \gamma^{-1} h \gamma \in \pi_1(M)$, and this cannot happen 
if $\pi_1(M)$ is
biorderable.  

(2) If $\gamma$ is represented by a regular fibre the 
result follows from part (1). 
Suppose then that $\gamma \in 
\pi_1(M)$ is a class
represented by an arbitrary fibre of the given Seifert structure. 
Evidently there is an integer
$\alpha > 0$ such that $\gamma^{\alpha}$ is represented by a regular 
fibre. Thus lemma \ref{central} shows that  $\gamma$ is central.
\qed

The proof of the ``necessity'' part of theorem \ref{SeifertBO} will 
now be divided into the three cases
$\Sigma = \emptyset, L \ne \emptyset$, and $\Sigma \ne L = \emptyset$.
  \\

\noindent {\bf Case 1:} $\Sigma = \emptyset$.

In this case $M \to B$ is an orientable, locally trivial circle bundle
(lemma \ref{orientable}). If $B \cong S^2$, then
$M$ is homeomorphic to either $S^3$, a lens space with finite, 
non-trivial fundamental group, or $S^1\times S^2$. Evidently the 
second option is incompatible with ($*$). Suppose then that $B$ is 
$P^2$ or $K$. Since $M$ satisfies ($*$), it is clear in these cases 
that $M \to B$ cannot be a trivial bundle, and this fact determines $M$ 
up to homeomorphism.
To see this we recall that the orientable circle bundles over
$B$ are classified  by the set of homotopy classes of maps $B \to BS^1$. Since $BS^1 = K(\mathbb Z, 2)$, these bundles correspond to elements in
$H^2(B) \cong \mathbb Z / 2$. In particular there is a unique, 
orientable, non-trivial circle bundle $p\co M \to B$. In order to 
construct $M$, let $D$ be a small $2$-disk in $B$ and set
$B_0 = \overline{B \setminus D}$. 
Consider  $M' = (B_0 \times S^1) \cup_f (D^2 \times S^1)$ where
$f\co \partial B_0 \times S^1 \to S^1 \times S^1$ preserves the $S^1$ 
factors and identifies $\partial D^2 \times pt$ with a curve in 
$\partial B_0 \times S^1$ which wraps once around $\partial B_0$ 
and once around $S^1$. There is a natural map $M' \to B$ which is an
orientable circle bundle over $B$ and it is simple to see that 
$H_1(M') \not \cong H_1(B \times S^1)$. Thus $M'$ is the bundle we are 
looking for: $M' \cong M$.

{\bf Subcase:}
$B = P^2$. The explicit description given in the previous paragraph
of the closed, connected, non-orientable manifold $M$ shows that $\pi_1(M)
\cong \Z$. Proposition \ref{sfsbasic} (5) shows that $M$ is one of the
manifolds in the given list, namely $S^1 \tilde \times S^2$.

{\bf Subcase:} $B = K$.  We will describe why this cannot happen. 
From  \S \ref{sfs} we see that
$$\pi_1(M) \cong \langle x,y,t \; | \; t \mbox{ central }, 
t = x^2y^2 \rangle  \; \cong \; 
\langle x,y \; | \; x^2y^2 \mbox{ central} \rangle.$$
We verify $x^2$ is central in this
group by the calculation
$$[x^2,y] = x^2yx^{-2}y^{-1} = (x^2y^2)y^{-1}x^{-2}y^{-1} =
y^{-1}x^{-2}(x^2y^2)y^{-1} = 1,$$
and so by lemma \ref{central}, $x$ is central as well. But this is easily 
seen to be false by projecting $\pi_1(M)$ onto the non-abelian group
$\langle  x,y \; | \; x^2, y^2 \rangle \; \cong 
\mathbb Z / 2 * \mathbb Z / 2$.  We've shown that if $M$ satisfies ($*$), 
it cannot be a circle bundle over the Klein bottle. \\


\noindent {\bf Case 2:} $L \ne \emptyset$, that is, there are reflector curves.

We will show that in this case, $M$ is either $S^1 \tilde \times S^2$, a
solid Klein bottle, or a trivial circle bundle over the M\"{o}bius band.

Let $N$ be a regular neighbourhood in $B$ of the set of reflector lines and
$N_0$ a component of $N$. Let $\gamma$ be the central element of
$\pi_1(M)$ represented by an exceptional fibre in $N_0$ (cf. lemma
\ref{orientable} (2)).
Set $B_0 = \overline{B \setminus N_0}$ and observe that the decomposition
$B = B_0 \cup N_0$ induces a splitting $M = M_0 \cup P_0$ where 
$M_0 \to B_0$ and $P_0 \to N_0$ are Seifert
fibrings. One readily verifies that $M_0 \cap P_0$ is a vertical torus or a
vertical annulus depending on whether $L \cap N_0$ is a circle or an arc (vertical Klein bottles are ruled out by lemma \ref{orientable}). 
It follows that $P_0$ is a twisted $I$-bundle over a torus in the first 
case and a solid Klein bottle (cf. pages 433-434
of \cite{Sc2}) otherwise. In any event, $M_0 \cap P_0$ is incompressible 
in $P_0$.

Now we distinguish two cases, namely the one where 
$M_0\cap P_0$ is compressible in $M_0$, and the one where it is not.
Since a fibre is never contractible in a Seifert manifold with boundary,
and the only Seifert manifolds with compressible boundaries are 
homeomorphic to solid tori or solid Klein bottles, our assumptions imply 
that if $M_0 \cap P_0$ compresses in $M_0$, then $M_0$ is a
solid torus, $P_0$ is a twisted $I$-bundle over the torus, and 
$M_0 \cap P_0 = \partial M_0$, i.e. $M$ is a Dehn filling of $P_0$. 
It follows that $\pi_1(M)$ is a non-trivial quotient group of 
$\pi_1(P_0) = \mathbb Z^2$. On the other hand, the bi-orderability of 
$\pi_1(M)$ implies it has no torsion.  The only possibility is for
$\pi_1(M) \cong \mathbb Z$.  Since $M$ is closed and non-orientable, it 
must be $S^1 \tilde \times S^2$.

Assume then that $M_0 \cap P_0$ is incompressible in
$M_0$, so that $\pi_1(M)$ is the free product of $\pi_1(M_0)$
and $\pi_1(P_0)$ amalgamated along
$\pi_1(M_0 \cap P_0)$. 
As $\gamma \in \pi_1(P_0) \setminus \pi_1(M_0 \cap P_0)$, the only way
it can be central is for $\pi_1(M_0 \cap P_0) \to \pi_1(M_0)$ to be an
isomorphism. It follows that $M \cong P_0$ and so is either a solid Klein 
bottle or twisted $I$-bundle over the torus, both of which have 
bi-orderable fundamental groups. Noting that the latter space is homeomorphic
to a trivial $S^1$-bundle over the M\"{o}bius band completes this part of
the proof of theorem \ref{SeifertBO}. \\


\noindent {\bf Case 3:} $\Sigma \ne L = \emptyset$, that is, 
there are isolated singular fibres.

Let the orders of the cone points in ${\cal B}$ be
$\alpha_1, \ldots  \alpha_n \geqslant 2$ ($n \geqslant 1$). We shall argue
that $M$ is homeomorphic to one of $S^3, S^1 \times S^2, 
S^1 \tilde \times S^2$, or the trivial bundle $D^2 \times S^1$.

\begin{lemma} \label{n<3} If $M$ satisfies ($*$) and the conditions of case 3 
hold,
then there are fewer than three cone points in ${\cal B}$.  \end{lemma}

\pf If $n \geqslant 3$, there are surjective homomorphisms
$$\pi_1(M) \to \pi_1^{orb}({\cal B}) \to \Delta(\alpha_1, \alpha_2,
\alpha_3)  = \;
\langle \bar \gamma_1, \bar \gamma_2, \bar \gamma_3 \; | \; \bar
\gamma_j^{\alpha_j} = 1, \bar \gamma_1
\bar \gamma_2 \bar \gamma_3 = 1 \rangle$$
where $\Delta(\alpha_1, \alpha_2, \alpha_3)$ is the $(\alpha_1, \alpha_2,
\alpha_3)$ triangle
group and $\bar \gamma_1, \bar \gamma_2, \bar \gamma_3$ are the images of
classes in $\pi_1(M)$
corresponding to the first three exceptional fibres.
Since $\bar \gamma_1$ and $\bar \gamma_2$ generate
the non-abelian group $\Delta(\alpha_1, \alpha_2, \alpha_3)$,
$\gamma_1$ is not central in $\pi_1(M)$, contradicting lemma 
\ref{orientable} (2). \qed

\begin{lemma}
If $M$ satisfies ($*$), and the conditions of case 3 hold, then the base
orbifold has $H_1(B)$ finite, so that $B$ is one of $D^2, S^2$ or $P^2$.
Moreover if $n=2$, then $B$ is $D^2$ or $S^2$.
\end{lemma}

\pf It is not hard to see that if either
$H_1(B)$ is infinite, or $n = 2$ and $B$ is non-orientable,
then there is a finite covering
$f\co \hat{B} \to {B}$ so that the pullback orbifold $\hat{\cal B}$ has at
least three cone points.  Let $M_f \to \hat{\cal B}$ be the pull-back of
$M \to {\cal B}$, via $f$, so that $M_f$ is a covering space of $M$, as
well as a Seifert fibre space over $\hat{\cal B}$.  If $\pi_1(M)$ is
bi-orderable, so
is $\pi_1(M_f)$, but that contradicts lemma \ref{n<3}. \qed

{\bf Subcase:} $B = P^2$ and $n = 1$.
Think of ${\cal B}$ as the union of a M\"{o}bius band without
singularities and a disk
containing exactly one cone point. From lemma \ref{orientable} it follows
that $M$ is a Dehn filling of the product of a M\"{o}bius band and $S^1$. 
But then, condition ($*$) implies that $\pi_1(M) \cong \mathbb Z$ and as 
$M$ is closed and non-orientable, it must be homeomorphic to 
$S^1 \tilde \times S^2$ (by proposition \ref{sfsbasic} (5)).

{\bf Subcase:} $B = S^2$ and $n = 1$ or $2$.  Then $M$ is the union of two
solid
tori, and the only such manifolds with bi-orderable groups are $S^3$ and $S^1
\times S^2$.

{\bf Subcase:} $B = D^2$ and $n = 1$ or $2$. When $n = 2$,
$\pi_1^{orb}({\cal B})
\cong \Z/{\alpha_1} * \Z/{\alpha_2}$ (as $L = \emptyset$) where
the class in
$\pi_1(M)$  represented by the first exceptional fibre projects to a generator
of $\Z/{\alpha_1}$ under the surjection $\pi_1(M) \to
\pi_1^{orb}({\cal B})$. As this class
is not central, this case does not arise (cf. lemma \ref{orientable} (2)). On
the other hand, if
$n = 1$ then $M \cong S^1 \times D^2$.  This completes the proof of the
present case and hence that of theorem \ref{SeifertBO}.
\qed


\section{Orderability and Sol manifolds} \label{osol}

The goal of this section is to investigate the orderability of the
fundamental groups of Sol manifolds, and in particular to prove 
theorem \ref{Solchar}. See pages 470--472 of \cite{Sc2} for 
background information of Sol manifolds. The interiors of such 
manifolds are covered by $\mathbb R^3$ and so are $P^2$-irreducible.

We recall from theorem 4.17 of \cite{Sc2} that every
compact, connected manifold $M$ whose interior admits a complete
Sol metric carries
the structure of a $2$-dimensional bundle over a $1$-dimensional orbifold
with a connected surface of non-negative Euler characteristic as generic
fibre. When
$\partial M \ne \emptyset$, this implies that $M$ is homeomorphic
to either a $3$-ball, a solid torus, a solid Klein bottle, the product of a
torus
with an interval, or a twisted $I$-bundle over the Klein bottle $K$. Theorem
\ref{Solchar} clearly holds in these cases,
so from now on we shall assume that $M$ is closed. Theorem 5.3 (i) of 
\cite{Sc2} yields several possibilities for the topology of $M$. However, 
one of them can be excluded -- the reader can verify that the union 
of two twisted $I$-bundles over the torus is double covered by $S^1 
\times S^1 \times S^1$ and so is not a Sol manifold. Thus denoting 
the torus by $T$ and the Klein bottle by $K$, we have that $M$ is either \\
\indent (i) a $T$- or $K$-bundle over the circle, or\\
\indent (ii) non-orientable and the union of two twisted $I$-bundles
over $K$, which are glued together along their Klein bottle boundaries, or\\
\indent (iii) orientable and the union of two twisted $I$-bundles
over $K$, which are glued together along their torus boundaries. \\
In cases (i) and (ii), $\pi_1(M)$ is LO by theorem \ref{homtoLO} and corollary
\ref{gencond}.

\begin{prop} Let $M$ be a closed, connected Sol manifold. \\
$(1)$ $\pi_1(M)$ is LO if and only if cases
$(i)$ or $(ii)$ arise, that is if and only if $M$ is either non-orientable
or orientable and a
torus bundle over the circle. \\
$(2)$ $\pi_1(M)$ is O if and only if $M$ is a torus
bundle over the circle whose monodromy in $GL_2(\mathbb Z)$ has at least
one positive eigenvalue.
\end{prop}

\pf (1) It remains to prove that an orientable manifold carrying the Sol
metric which is a union of two
twisted $I$-bundles over the Klein bottle cannot have an LO fundamental
group. Our proof is an adaptation of an idea of Bergman \cite{bergman2}.

The Klein bottle $K$ has fundamental group 
$\pi_1(K)=\langle m,l \ |$ $l^{-1}ml=m^{-1}\rangle$
(with $m$ and $l$ standing for meridian and longitude respectively); any
element in $\pi_1(K)$ can be written in the form $m^a l^b$ 
($a,b\in \mathbb{Z})$. We note that in any left-ordering of
$\pi_1(K)$ we have $m\ll |l|$, i.e. if $l^{\epsilon}>1$ for some
$\epsilon\in \{1,-1\}$, then
$m^n<l^\epsilon$ for all $n\in \mathbb{Z}$. (For if we had
$1<l^\epsilon<m^n$, it would follow that
$1>m^{-n} l^\epsilon = l^\epsilon m^n > 1\cdot 1=1$.) It follows that in
any left-ordering we have
$m\ll |m^a l^b|$ whenever $b\neq 0$. Observe that this condition
characterizes the subgroup of $\pi_1(K)$ generated by $m$.

Now we recall that our 3-manifold $M$ consists of two twisted $I$-bundles
$N_1, N_2$, and
$\pi_1(\partial N_i)\cong \mathbb{Z}^2$ is an index 2 subgroup of
$\pi_1(N_i)$ with generators
$l^2$ and $m$. With this choice of generators, the glueing map $f$ can be
described by an element
of $GL_2(\mathbb Z)$. Moreover, $\pi_1(M)$ is an amalgamated product
$\pi_1(N_1) *_f \pi_1(N_2)$.
Let's assume that this group is LO. By restriction, we obtain
left-orderings on $\pi_1(N_1)$ and
$\pi_1(N_2)$. In $\pi_1(M)$, the meridian $m_1\in \pi_1(N_1)$ is identified
with an element
$f(m_1)\in \pi_1(N_2)$. By the previous paragraph, $m_1 \ll |m_1^a
l_1^{2b}|$ for all $a, b\in \mathbb{Z}$ with $b\neq 0$ -- note that
$m_1^a l_1^{2b}$ lies in the boundary torus. Thus the same must be true
for $f(m_1)\in \pi_1(N_2)$, and it follows that $f(m_1)$ is a meridian of
$N_2$. In other words,
$f$ must glue meridian to meridian, and the $2\times 2$-matrix representing
$f$ is of the form $\pmatrix{1 & 0 \cr \ast & 1}$. It is well-known that
there are Seifert structures on $N_1$ and $N_2$ for which $m_1$ and $m_2$
are represented by
circle fibres in $\partial N_1$ and $\partial N_2$. Thus $M$ is Seifert
fibred, not Sol as hypothesized.

(2) There can be no $\pi_1$-injective Klein bottles in a manifold whose
group is O, so we are
reduced to the case of a torus bundle over the circle. Suppose that
$M$ is such a manifold with monodromy $A \in GL_2(\mathbb Z)$. There is 
an exact sequence
$$1 \to \mathbb Z^2 \to \pi_1(M) \to \mathbb Z \to 1$$
where the right hand $\mathbb Z$ acts on the left-hand $\mathbb Z^2$ by $A$.
Hence $\pi_1(M)$ is bi-orderable if and only if there is an bi-ordering on
$\mathbb Z^2$ whose positive cone $P$ is invariant under $A$.
If we think of $\mathbb Z^2$ as a subgroup of $\mathbb R^2$, then any
bi-ordering of $\mathbb
Z^2$ is defined by a line $L \subset \mathbb R^2$ through the origin; the
positive cone consists
of the elements of $\mathbb Z^2$ which lie in one of the components of
$\mathbb R^2 \setminus L$
as well as the elements of $\mathbb Z^2 \cap L$ which lie to one side of $0
\in L$. If one eigenvalue of $A$, say
$\lambda_1$, is positive, with an associated eigenvector $v_1 \in \mathbb
R^2$, there is a
linearly independent eigenvector $v_2 \in \mathbb R^2$ for $A$ whose
associated eigenvalue
$\lambda_2$ is real (since $\lambda_1 \lambda_2 = \pm 1$). We claim that
the positive
cone $P_L$ of the bi-order on $\mathbb Z^2$ defined by 
$L = \{ tv_2\;|\; t \in \mathbb R\}$ is invariant under the action of $A$. 
The fact that $M$ is Sol implies that the eigenvectors of
$A$ have irrational slopes - when $|A| = 1$ this follows from the fact that
$|\mbox{trace}(A)| >  2$, and when $|A| = -1$ from the fact that 
$|\mbox{trace}(A^2)|  > 2$. Hence $\mathbb Z^2 \cap L = \{0\}$ and so 
$P_L$ is the intersection of $\mathbb Z^2$ with a
component of $\mathbb R^2 \setminus L$. These components are preserved by
$A$ since $\lambda_1 >  0$, and thus $A(P_L) = P_L$.

On the other hand if the eigenvalues of $A$ are both negative,
then no half-space of $\mathbb R^2$ is preserved by $A$, and therefore
$\pi_1(M)$ admits no bi-ordering.
\qed

It follows from the description of the closed, connected Sol manifolds
we gave at the beginning of this section, that each such manifold is
finitely covered by a
torus bundle over the circle whose monodromy has positive eigenvalues. Thus,

\begin{cor}
The fundamental group of a closed, connected Sol manifold is virtually
bi-orderable.
\qed
\end{cor}


\section{Hyperbolic manifolds} \label{hyp}

Finally, we consider what is perhaps the most important 3-dimensional
geometry, and the least understood in terms of orderability.
R. Roberts, J. Shareshian, and M. Stein have very recently discovered a
family of closed hyperbolic 3-manifolds whose fundamental groups are 
not left-orderable.
These are constructed from certain fibre bundles over $S^1$, with fibre a
punctured torus, and pseudo-Anosov monodromy represented by the matrix
$\pmatrix{m & 1 \cr -1 & 0}$, where $m< -2$ is an odd negative integer.
The manifold $M^3_{p,q,m}$ is constructed by Dehn filling of this bundle,
corresponding to relatively prime integers $p > q \geqslant 1$.  We refer
the reader to \cite{RSS} for details of the construction.  In particular,
they show that
$$
\pi_1(M^3_{p,q,m}) \cong \langle t, a, b : t^{-1}at = aba^{m-1},
t^{-1}bt = a^{-1}, t^{-p} = (aba^{-1}b^{-1})^q \rangle,
$$ 
and prove that every homomorphism
$$
\pi_1(M^3_{p,q,m}) \to Homeo_+(\R)
$$ 
is trivial (in the sense defined in section \ref{cod1}).  
It follows that $\pi_1(M^3_{p,q,m})$ is not left-orderable. \\

\noindent {\bf Proof of Theorem \ref{geoms}}.  We need to show that each of
the eight
geometries contains manifolds whose groups are left-orderable and others whose
groups are not.   For the six Seifert geometries, this is an easy
consequence of theorem \ref{SeifertLO}.  First note that an
$S^3$-manifold has an LO group if and only if it is a $3$-sphere.  For
each of the other five Seifert geometries  one can construct prime,
orientable, closed manifolds with positive first
Betti number and which carry the appropriate geometric structure. Such
manifolds have LO groups by theorem \ref{homtoLO}.  On the other hand,
closed orientable manifolds admitting such geometries can
be constructed having first Betti number $0$ and non-orientable base
orbifold. Theorem \ref{SeifertLO} implies that their groups are not LO.
The case of closed manifolds admitting a Sol geometric structure can be
dealt with in a similar manner.  Likewise, there are many hyperbolic
closed manifolds with positive first Betti number, whose groups are
therefore LO.  Finally, the examples of \cite{RSS} provide many closed 
hyperbolic $3$-manifolds with non-LO groups. \qed


{\small
\vspace{5mm}
\noindent
D\'epartement de math\'ematiques, UQAM \\
P. O. Box 8888, Centre-ville \\
Montr\'eal, H3C 3P8 \\
Qu\'ebec, Canada \\
e-mail: boyer@math.uqam.ca

\vspace{3mm}
\noindent
Department of Mathematics, UBC \\
Room 121, 1984 Mathematics Road \\
Vancouver, V6T 1Z2 \\
B.C., Canada \\
e-mail: rolfsen@math.ubc.ca

\vspace{3mm}
\noindent
Institut Math\'ematique, Université de Rennes 1 \\
Campus de Beaulieu \\
35042 Rennes Cedex \\
France \\
e-mail: bertold.wiest@math.univ-rennes1.fr}


\begin{thebibliography}{FGRRW}

{\small

\bibitem[Ba]{Baumslag} G.\ Baumslag, {\it On generalised free products},
Math. Z. {\bf 78} (1962), 423--438

\bibitem[Be1]{bergman}
G.\ Bergman, {\it Left-orderable groups which are not locally indicable},
Pac.\ J.\ Math {\bf 147} (1991), 243--248

\bibitem[Be2]{bergman2}
G.\ Bergman, {\it Ordering coproducts of groups and semigroups},
J.\ Algebra {\bf 133} (1990), 313--339

\bibitem[BH]{BH}
R.\ Burns, V.\ Hale, {\it A note on group rings of certain torsion-free groups},
Can.\ Math.\ Bull.\ {\bf 15} (1972), 441--445.

\bibitem[Ca1]{Ca1}
D. Calegari, {\it $\R$-covered foliations of hyperbolic 3-manifolds}, Geometry and Topology {\bf 3} (1999), 137--153.

\bibitem[Ca2]{Ca2}
D. Calegari, {\it The geometry of $\mathbb R$-covered foliations},
Geom.\ Topol.\ {\bf 4} (2000), 457--515.

\bibitem[CJ]{CJ} 
A.\ Casson, D.\ Jungreis, {\it Convergence groups and Seifert fibered 
$3$-manifolds}, Invent.\ Math.\ {\bf 118} (1994), 441--456.

\bibitem[CD]{CD}
D.\ Calegari, N.\ Dunfield, {\it Laminations and groups of 
homeomorphisms of the circle},
preprint, 2002.

\bibitem[CR]{CR}
C.\ Champetier, V.\ Guirardel,
{\it Limit groups as limits of free groups:
compactifying the set of free groups}, ArXiv math. GR/0401042.

\bibitem[CK]{ChiswellKropholler}
I.\ Chiswell, P.\ Kropholler, {\it Soluble right orderable groups are 
locally indicable}, Canad.\ Math.\ Bull.
{\bf 36} (1993), 22--29.

\bibitem[Co]{Co}
P.\ F.\ Conrad, {\it Right-Ordered Groups}, Michigan Math.\ J., 
{\bf 6} (1959), 267--275.

\bibitem[De]{Dehornoy}
P.\ Dehornoy, {\it Braid groups and left distributive operations}, Trans.\
Amer.\ Math.\ Soc.\ {\bf 345} (1994), 115--150.

\bibitem[DT]{DT} N. Dunfield and W. Thurston, 
{\it The virtual Haken conjecture: experiments and examples}, Geom.\ 
Topol.\ {\bf 7} (2003), 399-441.

\bibitem[Du]{Dunwoody} M.\ Dunwoody, {\it An equivariant sphere theorem},
Bull.\ London Math.\ Soc.\ 17 (1985), 437--448

\bibitem[EHN]{EHN}
D.\ Eisenbud, U.\ Hirsch, W.\ Neumann, 
{\it Transverse foliations on Seifert bundles and self-homeomorphisms 
of the circle}, Comm.\ Math.\ Helv.\ {\bf 56} (1981), 638--660.

\bibitem[Ep]{Ep} D.\ B.\ A.\ Epstein, 
{\it Projective planes in 3-manifolds}, Proc.\ LMS {\bf 11} (1961), 
469--484

\bibitem[Fa]{Fa}
T.\ Farrell, {\it Right-orderable deck transformation groups}, Rocky 
Mtn.\ J.\ Math. {\bf 6} (1976), 441--447.

\bibitem[Fe]{Fe}
S.\ Fenley, {\it Anosov flows in $3$-maifolds}, Ann.\ Math.\ (2)
{\bf 139} (1992), 79--115.

\bibitem[FGRRW]{FGRRW}
R.\ Fenn, M.\ Greene, D.\ Rolfsen, C.\ Rourke, B.\ Wiest, {\it Ordering the
braid groups}, Pac.\ J.\ Math {\bf 191} 1999, 49--74.

\bibitem[Ga1]{Ga1}
D.\ Gabai, {\it Foliations and 3-manifolds}, Proc.\ Int.\ Cong.\ Math., 
Kyoto 1990, 609--619.

\bibitem[Ga2]{Ga2}
D.\ Gabai, {\it Convergence groups are Fuchsian groups}, Ann. of Math. 
{\bf 136} (1992), 447--510.

\bibitem[Ga3]{Ga}
D.\ Gabai, {\it Eight problems in the theory of foliations and laminations},
Geometric Topology, Volume $2$,
editor W. Kazez, AMS/IP Studies in Advanced Mathematics, 1996, 1--34.

\bibitem[GS]{GagSp} A.\ M.\ Gaglione, D.\ Spellman, {\it Generalisations
of free groups: some questions}, Comm. Algebra {\bf 22} 1994, 3159--3169.

\bibitem[G-M]{G-M}
J.\ Gonz\'alez-Meneses, {\it Ordering pure braid groups on closed surfaces},
Pac.\ J.\ Math {\bf 203} (2002), 369--378

\bibitem[GL]{GorinLin} E.\ A.\ Gorin and V.\ Ya.\ Lin,  
{\it Algebraic equations with continuous coefficients and some problems 
of the algebraic theory of braids}, Math. USSR Sbornik {\bf 7} (1969), 
569--596.

\bibitem[Ha]{Ha} A.\ Hatcher, {\it Notes on basic 3-manifold topology},
available electronically at http://www.math.cornell.edu/\~{}hatcher

\bibitem[He]{He}
J.\ Hempel, {$3$-manifolds}, Ann. of Math Studies {\bf 86}, 1976.

\bibitem[HiSt]{HiSt}
P.\ Hilton, U.\ Stammbach, {\it A course in homological algebra}, GTM 4,
Springer-Verlag, 1971.

\bibitem[HoSh]{HoSh}
J.\ Howie, H.\ Short, {\it The band-sum problem}, J.\ London Math.\ Soc.\ 
{\bf 31} (1985), 571-576.

\bibitem[Jc] {Jc}
W.\ Jaco, {\it Lectures on three-manifold topology},
CBMS Regional Conf.\ Ser.\ Math.\ {\bf 43} (1980).

\bibitem[Jn]{Jn} M. Jankins,
{\it The space of homomorphisms of a Fuchsian groups to 
$PSL_2(\mathbb R)$}, dissertation, University of Maryland, 1983.

\bibitem[JN1]{JN1}
M.\ Jankins, W.\ Neumann, 
{\it Homomorphisms of Fuchsian groups to $PSL_2(\mathbb R)$}, 
Comm.\ Math.\ Helv.\ {\bf 60} (1985), 480--495.

\bibitem[JN2]{JN2}
M.\ Jankins, W.\ Neumann,
{\it Rotation numbers and products of circle homomorphisms}, Math.\ Ann.\
{\bf 271} (1985), 381--400.

\bibitem[KR]{KR}
D.\ Kim, D.\ Rolfsen,
{\it  An ordering for groups of pure braids and fibre-type hyperplane 
arrangements},  Canad.\ J.\ Math.\ {\bf 55} (2003), 822--838.

\bibitem[LR]{LR}
R.\ H.\ La Grange, A.\ H.\ Rhemtulla, 
{\it A remark on the group rings of order preserving permutation groups,} 
Canad. Math. Bull. {\bf 11} (1968), 679-680.

\bibitem[Li]{Linnell}
P.\ Linnell, {\it Left ordered amenable and locally indicable groups}, J.
London Math. Soc. (2)
{\bf 60} (1999), 133--142.

\bibitem[Lo]{Long}
D.\ Long, {\it Planar kernels in surface groups},
Quart.\ J.\ Math.\ Oxford (2), {\bf 35} (1984), 305--310

\bibitem[Lu]{Lu}
J.\ Luecke,
{\it Finite covers of $3$-manifolds containing essential tori},
Trans.\ Amer.\ Math.\ Soc.\ {\bf 310} (1988), 381--391.

\bibitem[LS]{LySch}
R.\ Lyndon, P.\ Schupp, {\it Combinatorial group theory}. Ergebnisse der
Mathematik und ihrer Grenzgebiete, Band 89. Springer-Verlag,
Berlin-New York, 1977

\bibitem[MKS]{MKS}
W.\ Magnus, A.\ Karrass, and D.\ Solitar,
{\it Combinatorial group theory}, Dover, New York, revised edition, 1976.

\bibitem[Ma]{malcev}
A.\ I.\ Mal'cev, {\it On the embedding of group algebras in division 
algebras,} Dokl.\ Akad.\ Nauk SSSR {\bf 60} (1948) 1944-1501.

\bibitem[MSY]{MSY}
W.\ Meeks, L.\ Simon, and S.\ T.\ Yau,
{\it Embedded minimal surfaces, exotic spheres, and manifolds with 
positive Ricci curvature},  Ann.\ of Math. (2) {\bf 116} (1982), 621--659.

\bibitem[MR]{MR}
R.\ Botto Mura, A.\ H.\ Rhemtulla, {\it Orderable groups}, 
Lecture notes in pure and applied mathematics {\bf 27}, 
Marcel Dekker Inc., New York-Basel, 1977.

\bibitem[Na]{Na}
R. Naimi, {\it Foliations transverse to fibers of Seifert manifolds},
Comm.\ Math.\ Helv.\ {\bf 69} (1994), 155--162.

\bibitem[Ne]{neumann}
B.\ H.\ Neumann, {\it On ordered division rings,} 
Trans.\ Amer.\ Math.\ Soc.\ {\bf 66} (1949) 202-252.

\bibitem[No]{No}
S.\ P.\ Novikov, {\it Topology of foliations},
Trans.\ Moscow Math.\ Soc.\ {\bf 14} (1965), 268-304.

\bibitem[Pa]{Passman}
D.\ S.\ Passman, {\it The algebraic structure of group rings}, Pure and
applied mathematics, Wiley-Interscience, 1977.

\bibitem[PR]{PR} B.\ Perron and D.\ Rolfsen, 
{\it On orderability of fibred knot groups}, preprint.

\bibitem[RR]{RR}
A. Rhemtulla, D. Rolfsen, {\it Local indicability in ordered groups, braids
and elementary amenable groups}, preprint.

\bibitem[RSS]{RSS}
R.\ Roberts, J.\ Shareshian, M.\ Stein, personal communication.

\bibitem[RS]{RS}
R.\ Roberts, M.\ Stein, {\it Exceptional Seifert group actions on 
$\mathbb R$}, J.\ Knot Th.\ Ram.\ {\bf 8} (1999), 241--247.

\bibitem[RoWi]{RolW} D.\ Rolfsen, B.\ Wiest, {\it Free group automorphisms,
invariant orderings and applications}, Algebraic and Geometric Topology 1
(2001), 311--319.

\bibitem[RZ]{RZ}
D.\ Rolfsen, J.\ Zhu, {\it Braids, ordered groups and zero divisors},  
J.\ Knot Theory and Ramifications, {\bf 7} (1998), 837--841.

\bibitem[RW]{RW}
C.\ Rourke, B.\ Wiest, {\it Order automatic mapping class groups}, Pac.\ J.\
Math.\ {\bf 194}, (2000), 209--227.

\bibitem[Sc1]{Sc1}
G.\ P.\ Scott, {\it Compact submanifolds of $3$-manifolds}, J.\ London 
Math.\ Soc.\ {\bf 7} (1973), 246--250.

\bibitem[Sc2]{Sc2}
G.\ P.\ Scott,  {\it The geometries of $3$-manifolds}, 
Bull.\ Lond.\ Math.\ Soc.\ {\bf 15} (1983), 401--487.

\bibitem[Sc3]{Sc3} 
G.\ P.\ Scott, 
{\it There are no fake Seifert fibre spaces with infinite $\pi \sb{1}$}, 
Ann.\ of Math.\ {\bf 117} (1983), 35--70.

\bibitem[ScWa]{ScWa} G.\ P.\ Scott, T.\ Wall, {\it Topological methods in
group theory}, Homological group theory (Proc. Sympos., Durham, 1977),
137--203, London Math.\ Soc.\ Lecture Note Ser., 36,
Cambridge Univ. Press, Cambridge-New York, 1979.

\bibitem[Seif]{Seif} H.\ Seifert, 
{\it Topologie dreidimensionaler gefaserter R\"aume,} Acta Math.\ 
{\bf 60}(1933), 147--238.

\bibitem[SW]{SW}
H.\ Short, B.\ Wiest, {\it Orderings of mapping class groups after Thurston},
l'Ens.\ Math.\ 46 (2000), 279--312

\bibitem[Sm]{Smythe}
N.\ Smythe, {\it Trivial knots with arbitrary projection}, 
J.\ Austral.\ Math.\ Soc.\ {\bf 7} (1967), 481--489.

\bibitem[Th1]{Th1}
W.\ Thurston, {\it The geometry and topology of 3-manifolds}, Lecture notes,
Princeton University, 1977.

\bibitem[Th2]{Th2}
W.\ Thurston,  {\it Three-manifolds, foliations and circles, I}, preprint,
1997.

\bibitem[Tlf]{Tlf} J.\ Tollefson,  
{\it The compact $3$-manifolds covered by $S^2 \times \mathbb R$},
Proc.\ Amer.\ Math.\ Soc.\ {\bf 45} (1974), 461--462.
1997.

\bibitem[V]{vinogradov}
A.\ A.\ Vinogradov, {\it On the free product of ordered groups}, 
Mat.\ Sb.\ {\bf  67} (1949), 163--168.

}

\end{thebibliography}
\end{document}